\documentclass[fleqn]{mat01}
\usepackage{times,mathtimy,amssymb,latexsym,epsfig}
\begin{document}

\setcounter{page}{109} \firstpage{109}

\makeatletter
\def\artpath#1{\def\@artpath{#1}}
\makeatother \artpath{d:/prema/pm}

\font\xxxxx=msam10 at 10pt
\def\blacksquare{\mbox{\xxxxx{\char'245\ \!}}}

\newtheorem{theo}[defin]{\bf Theorem}
\newtheorem{lem}[defin]{Lemma}
\newtheorem{propo}[defin]{\rm PROPOSITION}

\def\remarr{\trivlist \item[\hskip \labelsep{\it Remark.}]}

\renewcommand\theequation{\thesection\arabic{equation}}

\markboth{P~K~Dutt, N~Kishore Kumar and
C~S~Upadhyay}{Nonconforming h-p spectral element methods}

\title{Nonconforming h-p spectral element methods for elliptic
problems}

\author{P~K~DUTT$^{*}$, N~KISHORE KUMAR$^{*}$ and C~S~UPADHYAY$^{\dagger}$}

\address{$^{*}$Department of Mathematics; $^{\dagger}$Department of Aerospace
Engineering, Indian Institute of Technology Kanpur,
Kanpur~208~016, India}

\volume{117}

\mon{February}

\parts{1}

\pubyear{2007}

\Date{MS received 20 February 2006; revised 23 July 2006}

\begin{abstract}
In this paper we show that we can use a modified version of the
h-p spectral element method proposed in
\cite{duttora1,duttom,duttora2,tomarth} to solve elliptic problems
with general boundary conditions to exponential accuracy on
polygonal domains using nonconforming spectral element functions.
A geometrical mesh is used in a neighbourhood of the corners. With
this mesh we seek a solution which minimizes the sum of a weighted
squared norm of the residuals in the partial differential equation
and the squared norm of the residuals in the boundary conditions
in fractional Sobolev spaces and enforce continuity by adding a
term which measures the jump in the function and its derivatives
at inter-element boundaries, in fractional Sobolev norms, to the
functional being minimized. In the neighbourhood of the corners,
modified polar coordinates are used and a global coordinate system
elsewhere. A~stability estimate is derived for the functional
which is minimized based on the regularity estimate in
\cite{babguo1}. We examine how to parallelize the method and show
that the set of common boundary values consists of the values of
the function at the corners of the polygonal domain. The method is
faster than that proposed in \cite{duttora1,duttom,tomarth} and
the h-p finite element method and stronger error estimates are
obtained.
\end{abstract}


\maketitle

\section{Introduction}

In \cite{duttora1,duttom,duttora2,tomarth} h-p spectral element
methods for solving elliptic boundary value problems on polygonal
domains using parallel computers were proposed. For problems with
Dirichlet boundary conditions the spectral element functions were
nonconforming. For problems with Neumann and mixed boundary
conditions the spectral element functions had to be continuous at
the vertices of the elements only. In this paper we propose a
modified version of this method using nonconforming spectral
element functions which works for general boundary conditions.

For simplicity of exposition we restrict ourselves to scalar
problems although the method applies to elliptic systems too.

A method for obtaining a numerical solution to exponential
accuracy for elliptic problems with analytic coefficients posed on
a curvilinear polygon whose boundary is piecewise analytic with
mixed Neumann and Dirichlet boundary conditions was first proposed
by Babuska and Guo \cite{babguo2} within the framework of the
finite element method. They were able to resolve the singularities
which arise at the corners by using a geometrical mesh. This
problem has also been examined by Karniadakis and Spencer in
\cite{karnia}.

We also use a geometrical mesh to solve the same class of problems
to exponential accuracy using h-p spectral element methods. In a
neighbourhood of the corners modified polar coordinates
$(\tau_{k},\theta_{k})$ are used, where $\tau_{k}=\ln\, r_{k}$ and
$(r_{k},\theta_{k})$ are polar coordinates with the origin at the
vertex $A_{k}$. Away from sectoral neighbourhoods of the corners a
global coordinate system is used consisting of $(x_{1},x_{2})$
coordinates.

We now seek a solution which minimizes the sum of the squares of a
weighted squared norm of the residuals in the partial differential
equation and the sum of the squares of the residuals in the
boundary conditions in fractional Sobolev norms and enforce
continuity by adding a term which measures the sum of the squares
of the jump in the function and its derivatives in fractional
Sobolev norms to the functional being minimized. These
computations are done using modified polar coordinates in sectoral
neighbourhoods of the corners and a global coordinate system
elsewhere in the domain. The spectral element functions are
nonconforming. For the modified version of the h-p spectral
element method examined here a stability estimate is proved which
is based on the regularity estimate of Babuska and Guo in
\cite{babguo1}. The proof is much simpler than that of the
stability estimate in \cite{duttora1,duttom}. Moreover the error
estimates are stronger.

The set of common boundary values for the numerical scheme
consists of the values of the function at the vertices of the
polygonal domain. Since the cardinality of the set of common
boundary values is so small we can compute a nearly exact
approximation to the Schur complement. Let $M$ denote the number
of corner layers and $W$ denote the number of degrees of freedom
in each independent variable of the spectral element functions,
which are a tensor product of polynomials, and let $W$ be
proportional to $M$. Then the method is faster than the h-p
spectral element method in \cite{duttora1,duttom,tomarth} by a
factor of $O(W^{1/2})$ and faster than the h-p finite element
method by a factor of $O(W)$.

We now outline the contents of this paper. In \S 2 function spaces
are defined and differentiability estimates are obtained. In \S 3
we state and prove stability estimates. In \S 4 the numerical
scheme, which is based on these estimates, is described and in \S
5 error estimates are obtained. In \S 6 we examine the issues of
parallelization and preconditioning. Finally \S 7 contains
technical results which are needed to prove the stability theorem.

\section{Function spaces and differentiability estimates}

Let $\Omega$ be a curvilinear polygon with vertices
$A_{1},A_{2},\dots,A_{p}$ and corresponding sides
$\Gamma_{1},\Gamma_{2},\dots,\Gamma_{p}$ where $\Gamma_{i}$ joins
the points $A_{i-1}$ and $A_{i}.$ We shall assume that the sides
$\overline{\Gamma}_{i}$ are analytic arcs, i.e.
\begin{equation*}
\overline{\Gamma}_{i}=\{
(\varphi_{i}(\xi),\psi_{i}(\xi))|\xi\in\overline{I}=[-1,1]\}
\end{equation*}
with $\varphi_{i}(\xi)$ and $\psi_{i}(\xi)$ being analytic
functions on $\overline{I}$ and
$|\varphi_{i}^{\prime}(\xi)|^{2}+|\psi_{i}^{\prime}(\xi)|^{2}\geq\alpha>0.$
By $\Gamma_{i}$ we mean the open arc, i.e. the image of
$I=(-1,1)$.

Let the angle subtended at $A_{j}$ be $\omega_{j}.$ We shall
denote the boundary $\partial\Omega$ of $\Omega$ by $\Gamma$.
Further, let $\Gamma=\Gamma^{[0]}\bigcup\Gamma^{[1]},$
$\Gamma^{[0]}=\bigcup_{i\in\mathcal{D}}\overline{\Gamma}_{i},$
$\Gamma^{[1]}=\bigcup_{i\in\mathcal{N}}\overline{\Gamma}_{i}$
where $\mathcal{D}$ is a subset of the set $\{ i| i=1,\dots,p\}$
and $\mathcal{N}=\{ i| i=1,\dots,p\} \setminus\mathcal{D}$. Let
$x$ denote the vector\break $x=(x_{1},x_{2})$.

Let $\mathcal{L}$ be a strongly elliptic operator
\begin{equation}
\frak{\mathcal{L}}(u)=-\sum_{r,s=1}^{2}(a_{r,s}(x)u_{x_{s}})_{x_{r}}+\sum_{r=1}^{2}b_{r}(x)u_{x_{r}}+c(x)u,\label{Eelloper}
\end{equation}
where $a_{s,r}(x)=a_{r,s}(x), b_{r}(x), c(x)$ are analytic
functions on $\overline{\Omega}$ and for any
$\xi_{1},\xi_{2}\in\Bbb{R}$ and any $x\in\overline{\Omega}$,
\begin{equation}
\sum_{r,s=1}^{2}a_{r,s}\xi_{r}\xi_{s}\geq\mu_{0}(\xi_{1}^{2}+\xi_{2}^{2})\label{Eellcond}
\end{equation}
with $\mu_{0}>0.$ In this paper we shall consider the boundary
value problem
\begin{align}
\mathcal{\frak{\mathcal{L}}}u & =
f\quad\mathrm{on} \,\ \Omega,\nonumber\\[.3pc]
u & = g^{[0]}\quad\mathrm{on} \,\ \Gamma^{[0]},\nonumber \\[.3pc]
\left(\frac{\partial u}{\partial N}\right)_{A} & =
g^{[1]}\quad\mathrm{on} \,\ \Gamma^{[1]},\label{Eellbvp}
\end{align}
where $({\partial u}/{\partial N})_{A}$ denotes the usual conormal
derivative which we shall now define. Let $A$ denote the
$2\times2$ matrix whose entries are given by
\begin{equation*}
A_{r,s}(x)=a_{r,s}(x)
\end{equation*}
for $r,s=1,2.$ Let $N=(N_{1},N_{2})$ denote the outward normal to
the curve $\Gamma_{i}$ for $i\in\mathcal{N}.$ Then
$\big(\frac{\partial u}{\partial N}\big)_{A}$ is defined as
follows:
\begin{equation}
\left(\frac{\partial u}{\partial
N}\right)_{A}\left(x\right)=\sum_{r,s=1}^{2}N_{r}a_{r,s}\frac{\partial
u}{\partial x_{s}}.
\end{equation}
Moreover let the bilinear form induced by the operator
$\mathcal{\frak{\mathcal{L}}}$ satisfy the inf--sup conditions. It
shall be assumed that the given data $f$ is analytic on
$\overline{\Omega}$ and $g^{[l]},l=0,1$ is analytic on every
closed arc $\overline{\Gamma}_{i}$ and $g^{[0]}$ is continuous on
$\Gamma^{[0]}$.

By $H^{m}(\Omega)$ we denote the Sobolev space of functions with
square integrable derivatives of order$\,\leq m$ on $\Omega$
furnished with the norm
\begin{equation*}
\Vert u\Vert_{_{H^{m}(\Omega)}}^{2}=\sum_{|\alpha|\leq m}\Vert
D^{\alpha}u\Vert_{_{L^{2}(\Omega)}}^{2}.
\end{equation*}

Define $r_{i}(x)$ to be the Euclidean distance between $x$ and the
vertex $A_{i}$ of $\Omega.$ Let
$\beta=(\beta_{1},\beta_{2},\dots,\beta_{p})$ denote a $p$-tuple
of real numbers, $0<\beta_{i}<1,$ $i=1,\dots,p.$ For any integer
$k$, let $\beta+k=(\beta_{1}+k,\beta_{2}+k,\dots,\beta_{p}+k).$
Further, we denote
\begin{equation*}
\Phi_{\beta}(x)=\prod_{i=1}^{p}r_{i}^{\beta_{i}}\quad\textrm{and}
\quad
\Phi_{\beta+k}(x)=\prod_{i=1}^{p}r_{i}^{\beta_{i}+k}.
\end{equation*}
Let $H_{_{\beta}}^{^{m,l}}(\Omega),m\geq l\geq0,$ $l$ an integer,
denote the completion of the set of all infinitely differentiable
functions under the norm
\begin{align*}
\Vert u\Vert_{_{H_{\beta}^{m,l}(\Omega)}}^{2} &=\Vert u\Vert
_{_{H^{l-1}(\Omega)}}^{2}+\sum_{|\alpha|=k,k=l}^{m}\Vert
D^{\alpha}u\,\Phi_{\beta+k-l}\Vert_{_{L^{2}(\Omega)}}^{2},\quad
l\geq1\\[.5pc]
\Vert u\Vert
_{_{H_{\beta}^{m,0}(\Omega)}}^{2}&=\sum_{|\alpha|=k,k=0}^{m}\Vert
D^{\alpha}u\,\Phi_{\beta+k-l}\Vert_{_{L^{2}(\Omega)}}^{2},\quad
l=0.
\end{align*}
For $m=l=0$ we shall write
$H_{\beta}^{0,0}(\Omega)=L_{\beta}(\Omega)$.\pagebreak

Let $\gamma$ be part of the boundary $\Gamma$ of $\Omega.$ Define
$H_{_{\beta}}^{^{m-\frac{1}{2},l-\frac{1}{2}}}(\gamma),\, m\geq
l,l\geq0$ to be the set of all functions $\phi$ on $\gamma$ such
that there exists $f\in H_{_{\beta}}^{^{m,l}}(\Omega)$ with
$\phi=f|_{\gamma}$ and
\begin{equation*}
\Vert \phi\Vert
_{_{H_{\beta}^{m-\frac{1}{2},l-\frac{1}{2}}(\gamma)}}=\inf_{_{f\in
H_{\beta}^{m,l}(\Omega)}}\{ \Vert f\Vert
_{_{H_{\beta}^{m,l}(\Omega)}}\}.
\end{equation*}
For $l$ an integer $0\leq l\leq2$, let
\begin{equation*}
\psi_{\beta}^{l}(\Omega)=\{ u(x)|u\in
H_{\beta}^{m,l}(\Omega),m\geq l\}
\end{equation*}
and
\begin{align*}
\frak{B}_{\beta}^{l}(\Omega) & = \{
u(x)|u\in\psi_{\beta}^{l}(\Omega),\Vert
\,|D^{\alpha}u|\Phi_{\beta+k-l}\Vert_{L^{2}(\Omega)}\leq
Cd^{k-l}(k-l)!\\[.3pc]
&\quad \ \mathrm{for}\:|\alpha|=k=l,l+1,\dots;d\geq1,\:
C\:\mathrm{independent}\:\mathrm{of}\: k\}.
\end{align*}

Let $Q\subseteq\Bbb{R}^{2}$ be an open set with a piecewise
analytic boundary $\partial Q$ and $\gamma$ be part or whole of
the boundary $\partial Q.$ Finally
$\frak{B}_{\beta}^{l-\frac{1}{2}}(\gamma),0\leq l\leq2,$ denotes
the space of all functions $\varphi$ for which there exists
$f\in\frak{B}_{\beta}^{l}(Q)$ such that $f=\varphi$ on $\gamma$.

Next as in \cite{babguo2} we introduce the space
$\frak{C}_{\beta}^{2}$:
\begin{align*}
\frak{C}_{\beta}^{2}(\Omega) & = \{ u\in
H_{\beta}^{2,2}(\Omega)|\,|D^{\alpha}u(x)|\leq
Cd^{k}k!(\Phi_{k+\beta-1}(x))^{-1},\\[.3pc]
&\quad \
|\alpha|=k=1,2,\dots;C\geq1;d\geq1\:\mathrm{independent}\:\mathrm{of}\:
k\} .
\end{align*}
The relationship between $\frak{C}_{\beta}^{2}$ and
$\frak{B}_{\beta}^{2}$ is given by Theorem 2.2 of \cite{babguo2}
which can be stated as follows:
\begin{equation*}
\frak{B}_{\beta}^{2}(\Omega)\subseteq\frak{C}_{\beta}^{2}(\Omega).
\end{equation*}

We need to state our regularity estimates in terms of local
variables which are defined on a geometrical mesh imposed on
$\Omega$ as in \S5 of \cite{babguo2}. $\Omega$ is first divided
into subdomains. Thus we divide $\Omega$ into $p$ subdomains
$S^{1},\dots,S^{p},$ where $S^{i}$ denotes a domain which contains
the vertex $A_{i}$ and no other, and on each $S^{i}$ we define a
geometrical mesh. Let $\frak{S}^{k}=\{ \Omega_{i,j}^{k},j=
1,\dots,J_{k},i=1,\dots,I_{k,j}\} $ be a partition of $S^{k}$ and
let $\frak{S}=\bigcup_{k=1}^{p}\frak{S}^{k}.$ The geometrical mesh
imposed on $\Omega$ is as shown in figure~1.

\begin{figure}[t]
\hskip 4pc{\epsfxsize=5.5cm\epsfbox{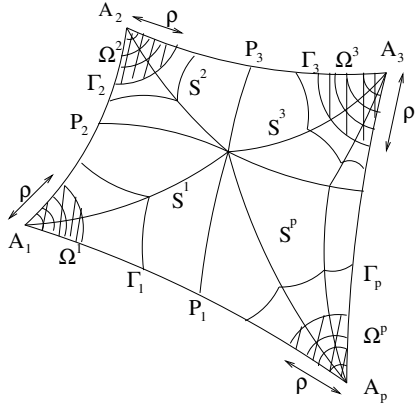}}\vspace{-.7pc}
\caption{Geometrical mesh with $M$ layers in the radial
direction.}\vspace{.7pc}
\end{figure}

We now put some restrictions on $\frak{S}.$ Let
$(r_{k},\theta_{k})$ denote polar coordinates with center at
$A_{k}.$ Let $\tau_{k}=\ln r_{k}.$ Choose $\rho$ so that the
curvilinear sector $\Omega^{k}$ with sides $\Gamma_{k}$ and
$\Gamma_{k+1}$ bounded by the circular arc $B_{\rho}^{k},$ center
at $A_{k}$ and radius $\rho$ satisfies
\begin{equation*}
\Omega^{k}\subseteq\bigcup_{\Omega_{i,j}^{k}\in\frak{S}^{k}}\overline{\Omega}_{i,j}^{k}.
\end{equation*}
$\Omega^{k}$ may be represented as
\begin{equation}
\Omega^{k}=\{ (x_{1},x_{2})\in\Omega\hbox{\rm :}\ 0<r_{k}<\rho\
\}.
\end{equation}

Let $\gamma_{i,j,l}^{k}, 1\leq l\leq4$ be the side of the
quadrilateral $\Omega_{i,j}^{k}\in\frak{S}.$ Then it is assumed
that
\begin{subequations}
\begin{align}
&\gamma_{i,j,l}^{k}\hbox{\rm :}\ \begin{cases}
x_{1} = h_{i,j}^{k}\varphi_{i,j,l}^{k}(\xi),\\[.5pc]
x_{2} = h_{i,j}^{k}\psi_{i,j,l}^{k}(\xi),\end{cases}
-1\leq\xi\leq1, \ \ l=1,3\\[.5pc]
&\gamma_{i,j,l}^{k}\hbox{\rm :}\ \begin{cases}
x_{1} = h_{i,j}^{k}\varphi_{i,j,l}^{k}(\eta),\\[.5pc]
x_{2} =
h_{i,j}^{k}\psi_{i,j,l}^{k}(\eta),\end{cases}-1\leq\eta\leq1, \ \
l=2,4
\end{align}
\end{subequations}
and that for some $C\geq1$ and $L\geq1$ independent of $i,j,k$ and
$l$,
\begin{equation}
\left|\frac{d^{t}}{ds^{t}}\varphi_{i,j,l}^{k}(s)\right|,
\left|\frac{d^{t}}{ds^{t}}\psi_{i,j,l}^{k}(s) \right| \leq
CL^{t}t!, \quad t=1,2,\dots.
\end{equation}
We shall place further restrictions on the geometric mesh imposed
on $\Omega^{k}$ later. Some of the elements may be curvilinear
triangles.

\begin{figure}[b]\vspace{.7pc} 
\centerline{\epsfxsize=12.3cm\epsfbox{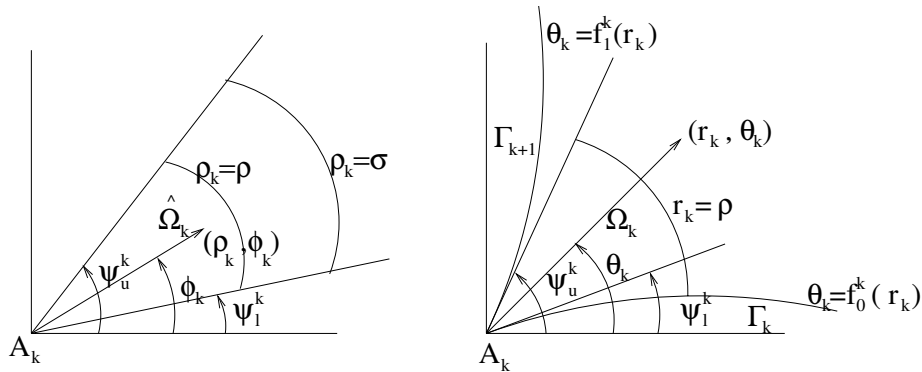}}\vspace{-.7pc}
\caption{Curvilinear sectors.}
\end{figure}

Let $(r_{k},\theta_{k})$ be polar coordinates with center at
$A_{k}.$ Then $\Omega^{k}$ is the open set bounded by the
curvilinear arcs $\Gamma_{k},$ $\Gamma_{k+1}$ and a portion of the
circle $r_{k}=\rho.$ We divide $\Omega^{k}$ into curvilinear
rectangles by drawing $M$ circular arcs
$r_{k}=\sigma_{j}^{k}=\rho\mu_{k}^{M+1-j},j=2,\dots,M+1,$ where
$\mu_{k}<1$ and $I_{k}-1$ analytic curves $C_{2},\dots,C_{I_{k}}$
whose exact form shall be prescribed in what follows. Let
$\sigma_{1}^{k}=0.$ Thus $I_{k,j}=I_{k}$ for $j\leq M;$ in fact,
we shall let $I_{k,j}=I_{k}$ for $j\leq M+1.$ Moreover
$I_{k,j}\leq I$ for all $k,j$ where $I$ is a fixed constant. Let
\begin{equation*}
\Gamma_{k+j}=\{
(r_{k},\theta_{k})|\theta_{k}=f_{j}^{k}(r_{k}),\quad
0<r_{k}<\rho\},
\end{equation*}
$j=0,1$ in a neighbourhood of $A_{k}$ in $\Omega^{k}.$ Then the
mapping
\begin{equation}
r_{k}=\rho_{k},\theta_{k}=\frac{1}{(\psi_{u}^{k}-\psi_{l}^{k})}[(\phi_{k}-\psi_{l}^{k})f_{1}^{k}(\rho_{k})-(\phi_{k}-\psi_{u}^{k})f_{0}^{k}(\rho_{k})],\label{Egamdef}
\end{equation}
where $f_{j}^{k}$ is analytic in $\rho_{k}$ for $j=0,1$, maps
locally the cone
\begin{equation*}
\{ (\rho_{k},\phi_{k})\hbox{\rm :}\
0<\rho_{k}<\sigma,\psi_{l}^{k}<\phi_{k}<\psi_{u}^{k}\}
\end{equation*}
onto a set containing $\Omega^{k}$ as in \S 3 of \cite{babguo2}.
The functions $f_{j}^{k}$ satisfy$f_{0}^{k}(0)=\psi_{l}^{k},$
$f_{1}^{k}(0)=\psi_{u}^{k}$ and $(f_{j}^{k})^{\prime}(0)=0$ for
$j=0,1.$ It is easy to see that the mapping defined in
(\ref{Egamdef}) has two bounded derivatives in a neighbourhood of
the origin which contains the closure of the open set
\begin{equation*}
\hat{\Omega}^{k}=\{ (\rho_{k},\phi_{k})\hbox{\rm :}\
0<\rho_{k}<\rho,\psi_{l}^{k}<\phi_{k}<\psi_{u}^{k}\}.
\end{equation*}
We choose the $I_{k}-1$ curves $C_{2},\dots,C_{I_{k}}$ as
\begin{equation*}
C_{i}\hbox{\rm :}\ \phi_{k}(r_{k},\theta_{k})=\psi_{i}^{k}
\end{equation*}
for $i=2,\dots,I_{k}.$ Here
\begin{equation*}
\psi_{l}^{k}=\psi_{1}^{k}<\psi_{2}^{k}<\cdots<\psi_{I_{k}+1}^{k}=\psi_{u}^{k}.
\end{equation*}
Let $\Delta\psi_{i}^{k}=\psi_{i+1}^{k}-\psi_{i}^{k}.$ Then $\{
\psi_{i}^{k}\}_{i,k}$ are chosen so that
\begin{equation}
\max_{i,k}(\Delta\psi_{i}^{k})<\lambda(\min_{i,k}(\Delta\psi_{i}^{k}))
\end{equation}
for some constant $\lambda.$ Another set of local variables
$(\tau_{k},\theta_{k})$ is needed in a neighbourhood of
$\Omega^{k}$ where
\begin{equation*}
\tau_{k}=\ln r_{k}.
\end{equation*}
In addition, we need one final set of local variables
$(\nu_{k},\phi_{k})$ in the cone
\begin{equation*}
\{ (\rho_{k},\phi_{k})\hbox{\rm :}\
0\leq\rho_{k}\leq\rho,\psi_{l}^{k}\leq\phi_{k}\leq\psi_{u}^{k}\},
\end{equation*}
where
\begin{equation*}
\nu_{k}=\ln\rho_{k}.
\end{equation*}
Let $S_{\mu}^{k}=\{ (r_{k},\theta_{k})\hbox{\rm :}\ 0\leq
r_{k}\leq\mu\} \cap\Omega.$ Then the image $\hat{S}_{\mu}^{k}$ in
$(\nu_{k},\phi_{k})$ variables of $S_{\mu}^{k}$ is given by
\begin{equation*}
\hat{S}_{\mu}^{k}=\{ (\nu_{k},\phi_{k})\hbox{\rm :}\
-\infty\leq\nu_{k}\leq\ln\mu,\psi_{l}^{k}\leq\phi_{k}\leq\psi_{u}^{k}\}.
\end{equation*}
Now the relationship between the variables $(\tau_{k},\theta_{k})$
and $(\nu_{k},\phi_{k})$ is given by
$(\tau_{k},\theta_{k})=M^{k}(\nu_{k},\phi_{k})$, viz.
\begin{align}
\tau_{k} & = \nu_{k},\nonumber\\[.1pc]
\theta_{k} & =
\frac{1}{(\psi_{u}^{k}-\psi_{l}^{k})}[(\phi_{k}-\psi_{l}^{k})f_{1}^{k}(\hbox{e}^{\nu_{k}})-(\phi_{k}-\psi_{u}^{k})f_{0}^{k}(\hbox{e}^{\nu_{k}})].
\end{align}
Hence it is easy to see that $J^{k}(\nu_{k},\phi_{k})$, the
Jacobian of the above transformation, satisfies $C_{1}\leq$
$|J^{k}(\nu_{k},\phi_{k})|\leq C_{2}$ for all
$(\nu_{k},\phi_{k})\in\hat{S}_{\mu}^{k}$, for all $0<\mu\leq\rho$.

We now need the fundamental regularity result from \cite{babguo1},
viz. Theorem 2.1 which we state as follows:

If $f\in H_{\beta}^{m,0}(\Omega)$, $g^{[j]}\in
H_{\beta}^{m+\frac{3}{2}-j,\frac{3}{2}-j}(\Gamma^{[j]})$, $j=0,1$,
$0<\beta_{i}<1$, $\beta_{i}>\beta_{i}^{^{\star}}$ and $m\geq0$,
then the solution of (2.3) exists in $H_{\beta}^{m+2,2}(\Omega)$
and
\begin{equation*}
\Vert u\Vert_{_{H_{\beta}^{m+2,2}(\Omega)}}\leq C_{m} \left( \Vert
f\Vert _{_{H_{\beta}^{m,0}(\Omega)}}+\sum_{j=0}^{1}\Vert
g^{[j]}\Vert_{_{H_{\beta}^{m+\frac{3}{2}-j,\frac{3}{2}-j}(\Gamma^{[j]})}}\right).
\end{equation*}
Let us define $\alpha_{i}=1-\beta_{i}^{^{\star}}$.

We now state the differentiability estimates for the solution $u$
of (\ref{Eellbvp}) which will be needed in this paper.

\begin{propo}$\left.\right.$\vspace{.5pc}

\noindent Let $1-\alpha_{k}>0.$ Then for
$\lambda_{k}<\alpha_{k}${\rm ,}
\begin{align}
&\int_{\psi_{l}^{k}}^{\psi_{u}^{k}}\int_{-\infty}^{\ln\mu}\sum_{|\varepsilon|\leq
m}|D_{\nu_{k}}^{\varepsilon_{1}}D_{\phi_{k}}^{\varepsilon_{2}}(u-u(A_{k}))|^{2}{\rm
e}^{-2\lambda_{k}\nu_{k}}\, {\rm d}\nu_{k}{\rm
d}\phi_{k}\nonumber\\[.4pc]
&\quad\, \leq\mu^{2\gamma_{k}}\,(C\,
d^{m-2}(m-2)!)^{2}\label{2.11}
\end{align}
for $0<\mu\leq\rho$ with $\gamma_{k}<\alpha_{k}-\lambda_{k}.$ If
$1-\alpha_{k}<0$ then for $\lambda_{k}<1/2${\rm ,} $(2.11)$
remains valid for $0<\mu\leq\rho$ with $\gamma_{k}=1/2.$
\end{propo}

The proposition can be proved in the same way as Theorem~{\rm 2.1}
of {\rm \cite{duttora1}}.\hfill $\blacksquare$

\section{The stability estimate}

\setcounter{equation}{0}

\setcounter{defin}{0}

Let
\begin{equation}
\mathcal{L}(u)=-\sum_{i,j=1}^{2}(a_{i,j}u_{x_{j}})_{x_{i}}+\sum_{i=1}^{2}b_{i}u_{x_{i}}+cu\label{1.1}
\end{equation}
be a strongly elliptic operator. We now consider the following
mixed boundary value problem:
\begin{align}
\mathcal{L}u &=f\quad\textrm{in} \ \ \Omega,\nonumber\\[.3pc]
\overline{\gamma}_{0}u&=u|_{\Gamma^{[0]}}=g^{[0]}\quad\textrm{and}\nonumber\\[.3pc]
\overline{\gamma}_{1}u&=\left(\frac{\partial u}{\partial
N}\right)_{A}\bigg|_{\Gamma^{[1]}}=g^{[1]}.\label{1.2}
\end{align}
Here the conormal derivative $\bar{\gamma}_{1}u$ is defined as
follows. Let $\gamma_{i}\subseteq\Gamma^{[1]}$ and let
$N=(N_{1},N_{2})^{T}$ denote the unit outward normal at a point on
$\gamma_{i}$. Then
\begin{equation}
\bar{\gamma}_{1}u=\left(\frac{\partial u}{\partial
N}\right)_{A}=\sum_{i,j=1}^{2}N_{i}\, a_{i,j}\,
u_{x_{j}}.\label{1.3}
\end{equation}
Moreover, let the bilinear form induced by the operator
$\mathcal{\frak{\mathcal{L}}}$ satisfy the inf--sup conditions.

We can now state the regularity result Theorem 2.1 of
\cite{babguo1} as follows:

Let $u$ be the solution to (3.2). Then
\begin{equation}
\Vert u\Vert_{_{H_{\beta}^{k+2,2}(\Omega)}}\leq C_{k}\left(\Vert
f\Vert_{_{H_{\beta}^{k,0}(\Omega)}}+\sum_{j=0}^{1}\Vert
g^{[j]}\Vert
_{_{H_{\beta}^{k+\frac{3}{2}-j,\frac{3}{2}-j}(\Gamma^{[j]})}}\right).\label{3.4}
\end{equation}
The above estimate for $k=0$ is used to prove the stability
estimate Theorem~3.1.

We remark that in Theorem~5.2 of \cite{babguo4}, Guo and Babuska
have extended the above regularity result to elliptic systems.
Hence the method applies to elliptic systems too.

Divide the polygonal domain $\Omega$ into $p$ sectors
$\Omega^{1},\Omega^{2},\dots,\Omega^{p}$ and a remaining portion
$\Omega^{p+1}.$ Further divide each of these subdomains into still
smaller elements
\begin{equation*}
\{ \Omega_{i,j}^{k},\;1\leq i\leq I_{k,j},\;1\leq j\leq M,\;1\leq
k\leq p\}.
\end{equation*}
Let
\begin{equation*}
\Omega^{p+1}=\{ \Omega_{i,j}^{k}\hbox{\rm :}\ 1\leq k\leq
p,M<j\leq J_{k},1\leq i\leq I_{k,j}\}.
\end{equation*}
We shall relabel the elements of $\Omega^{p+1}$ and write
\begin{equation*}
\Omega^{p+1}=\{ \Omega_{l}^{p+1}\hbox{\rm :}\ 1\leq l\leq L\}.
\end{equation*}
Now define the space of spectral element functions $\Pi^{M,W}=\{
\{ u_{i,j}^{k}(\nu_{k},\phi_{k})\}_{i,j,k},$ $\{
u_{l}^{_{^{p+1}}}(\xi,\eta)\}_{l}\},$ where $u_{i,1}^{k}=h_{k}$ a
constant for all $i$ and
\begin{equation*}
u_{i,j}^{k}(\nu_{k},\phi_{k})=\sum_{r=1}^{W_{j}}\sum_{s=1}^{W_{j}}g_{r,s}\,\nu_{k}^{r}\,\phi_{k}^{s},\quad
1<j\leq M.
\end{equation*}
Here $1\leq W_{j}\leq W.$ Moreover there is an analytic mapping
$M_{l}^{p+1}$ from the master square $S=(-1,1)^{2}$ to
$\Omega_{l}^{p+1}.$ We define
\begin{equation*}
u_{l}^{p+1}(M_{l}^{p+1}(\xi,\eta))=\sum_{r=1}^{W}\sum_{s=1}^{W}g_{r,s}\,\xi^{r}\,\eta^{s}.
\end{equation*}
Let $w\in H_{_{\beta}}^{^{2,2}}(\Omega)$. Now for $1\leq j\leq
M$,
\begin{equation}
\int_{\Omega_{i,j}^{k}}r_{_{k}}^{2\beta_{k}}|\mathcal{L}w|^{2}\hbox{d}x=\int_{\tilde{\Omega}_{i,j}^{k}}r_{_{k}}^{2(-1+\beta_{k})}|\tilde{\mathcal{L}}^{k}w|^{2}\,
\hbox{d}\tau_{k}\hbox{d}\theta_{k}.\label{1.15}
\end{equation}
Here $\tilde{\Omega}_{i,j}^{k}$ is the image of $\Omega_{i,j}^{k}$
in $(\tau_{k},\theta_{k})$ coordinates and
$\tilde{\mathcal{L}}^{k}w=r_{_{k}}^{2}\mathcal{L}w$. It has been
shown in \cite{duttom} that if we let $y_{1}=\tau_{k}$ and
$y_{2}=\theta_{k}$ then
\begin{equation}
\tilde{\mathcal{L}}^{k}w=-\sum_{i,j=1}^{2}\frac{\partial}{\partial
y_{i}}\left(\tilde{a}_{i,j}^{k}\,\frac{\partial w}{\partial
y_{j}}\right)+\sum_{i=1}^{2}\tilde{b}_{i}^{k}\,
w_{y_{i}}+\tilde{c}^{k}w.\label{1.16}
\end{equation}
Let $O^{k}$ denote the matrix
\begin{equation*}
O^{k}=\begin{bmatrix}
\cos\theta_{k} & -\sin\theta_{k}\\[.2pc]
\sin\theta_{k} & \cos\theta_{k}\end{bmatrix}
\end{equation*}
and
\begin{equation*}
\tilde{A}^{k}=\begin{bmatrix}
\tilde{a}_{1,1}^{k} & \tilde{a}_{1,2}^{k}\\[.5pc]
\tilde{a}_{2,1}^{k} & \tilde{a}_{2,2}^{k}\end{bmatrix}.
\end{equation*}
Then $\tilde{A}^{k}=(O^{k})^{T}AO^{k}$.

Let $J^{k}(\nu_{_{k}},\phi_{_{k}})$ denote the Jacobian of the map
$M^{k}(\nu_{k},\phi_{k})$ defined in \S 2. Then for $1<j\leq
M$,
\begin{equation}
\int_{\Omega_{i,j}^{k}}r_{_{k}}^{2\beta_{k}}|\mathcal{L}w|^{2}\hbox{d}x=\int_{\hat{\Omega}_{i,j}^{k}}\hbox{e}^{-2(1-\beta_{k})\nu_{k}}|\mathcal{L}_{i,j}^{k}w(\nu_{_{k}},\phi_{_{k}})|^{2}\,
\hbox{d}\nu_{k}\hbox{d}\phi_{k}.\label{1.7}
\end{equation}
Here $\hat{\Omega}_{i,j}^{k}$ is the image of $\Omega_{i,j}^{k}$
in $(\nu_{k},\phi_{k})$ variables and
\begin{equation*}
\mathcal{L}_{i,j}^{k}w=\sqrt{J^{k}}\,\,\tilde{\mathcal{L}}^{k}w.
\end{equation*}
Now
\begin{align*}
\mathcal{L}_{i,j}^{k}w(\nu_{k},\phi_{k}) &=
A_{i,j}^{k}w_{\nu_{k}\nu_{k}}+2B_{i,j}^{k}w_{\nu_{k}\phi_{k}}+C_{i,j}^{k}w_{\phi_{k}\phi_{k}}\\[.3pc]
&\quad\,
+D_{i,j}^{k}w_{\nu_{k}}+E_{i,j}^{k}w_{\phi_{k}}+F_{i,j}^{k}w.
\end{align*}

Let $\hat{A}_{i,j}^{k}$ be the polynomial approximation of
$A_{i,j}^{k},$ of degree $W_{j}$ in $\nu_{k}$ and $\phi_{k}$
separately, as defined in Theorem 4.46 of \cite{schwab}. Now we
define a differential operator with polynomial coefficients
$(\mathcal{L}_{i,j}^{k})^{^{a}},$ which is an approximation to
$\mathcal{L}_{i,j}^{k}$ as follows:
\begin{align*}
(\mathcal{L}_{i,j}^{k})^{a}w &=
\hat{A}_{i,j}^{k}w_{\nu_{k}\nu_{k}}+2\hat{B}_{i,j}^{k}w_{\nu_{k}\phi_{k}}+\hat{C}_{i,j}^{k}w_{\phi_{k}\phi_{k}}+\hat{D}_{i,j}^{k}w_{\nu_{k}}\\[.4pc]
&\quad\, +\hat{E}_{i,j}^{k}w_{\phi_{k}}+\hat{F}_{i,j}^{k}w.
\end{align*}
Let $\lambda_{k}=1-\beta_{k}.$ Then for $1<j\leq
M$,
\begin{subequations}
\begin{align}
&\left|\int_{\hat{\Omega}_{i,j}^{k}}|\mathcal{L}_{i,j}^{k}\,
w(\nu_{k},\phi_{k})|^{2}\hbox{e}^{-2\lambda_{k}\nu_{k}}\,
\hbox{d}\nu_{k}\hbox{d}\phi_{k}\right.\nonumber\\[.4pc]
&\qquad\, \left.
-\int_{\hat{\Omega}_{i,j}^{k}}|(\mathcal{L}_{i,j}^{k})^{^{a}}w(\nu_{k},\phi_{k})|^{2}\hbox{e}^{-2\lambda_{k}\nu_{k}}\,
\hbox{d}\nu_{k}\hbox{d}\phi_{k}\right|\nonumber\\[.5pc]
 &\quad\, \leq \varepsilon_{_{W}}((\rho\mu_{k}^{M+1-j})^{-2\lambda_{k}}\Vert w(\nu_{k},\phi_{k})-w(A_{k})\Vert_{_{2,\hat{\Omega}_{i,j}^{k}}}^{2}\nonumber\\[.5pc]
 &\qquad\,  +(\rho\mu_{k}^{M+1-j})^{4-2\lambda_{k}}|w(A_{k})|^{2}).\label{eq:1b}
 \end{align}
Here $\varepsilon_{_{W}}\rightarrow0$ as $W$$\rightarrow$$\infty$
and, in fact, $\varepsilon_{_{W}}$ is exponentially small in $W.$

Moreover, if $w(\nu_{k},\phi_{k})=w(A_{k}),$ a constant in
$\Omega_{i,1}^{k}$ for $1\leq i\leq I_{k},$ then
\begin{equation*}
\sum_{i=1}^{I_{k}}\int_{\hat{\Omega}_{i,1}^{k}}|\mathcal{L}_{i,1}^{k}w(\nu_{k},\phi_{k})|^{2}\hbox{e}^{-2\lambda_{k}\nu_{k}}\,
\hbox{d}\nu_{k}\hbox{d}\phi_{k}\leq\varepsilon_{_{M}}|w(A_{k})|^{2}.
\end{equation*}
Here $\varepsilon_{_{M}}\rightarrow0$ as $M\rightarrow\infty$ and
$\varepsilon_{_{M}}$ is exponentially small in $M$.

Hence we conclude that if $w(\nu_{k},\phi_{k})=w(A_{k}),$ a
constant in $\Omega_{i,1}^{k}$ for $1\leq i\leq I_{k},$
then
\begin{align}
& \sum_{i=1}^{I_{k}}\sum_{j=1}^{M}\int_{\hat{\Omega}_{i,j}^{k}}|\mathcal{L}_{i,j}^{k}w_{i,j}^{k}(\nu_{k},\phi_{k})|^{2}{\rm e}^{-2\lambda_{k}\nu_{k}}\, {\rm d}\nu_{k}{\rm
d}\phi_{k}\nonumber\\[.4pc]
&\quad\, \leq
C\,\left(\sum_{i=1}^{I_{k}}\sum_{j=2}^{M}(\rho\mu_{k}^{M+1-j})^{-2\lambda_{k}}\left(\int_{\hat{\Omega}_{i,j}^{k}}|(\mathcal{L}_{i,j}^{k})^{^{a}}w_{i,j}^{k}(\nu_{k},\phi_{k})|^{2}\,
{\rm d}\nu_{k}{\rm d}\phi_{k}\right)\right)\nonumber\\[.4pc]
 &\qquad\, +
\varepsilon_{_{W}}\,\left(\sum_{i=1}^{I_{k}}\sum_{j=2}^{M}(\rho\mu_{k}^{M+1-j})^{-2\lambda_{k}}\Vert
w_{i,j}^{k}-w(A_{k})\Vert
_{_{2,\hat{\Omega}_{i,j}^{k}}}^{2}+\left|w(A_{k})\right|^{2}\right)\nonumber\\[.4pc]
&\qquad\,
+\varepsilon_{_{M}}\left|w(A_{k})\right|^{2}.\label{0.1b}
\end{align}
\end{subequations}
Here $C$ is a constant.

Now
\begin{equation*}
\int_{\Omega_{l}^{p+1}}|\mathcal{L}w|^{2}\hbox{d}x_{1}\hbox{d}x_{2}=\int_{S}|\mathcal{L}w_{l}^{p+1}|^{2}J_{l}^{p+1}\,
\hbox{d}\xi \hbox{d}\eta.
\end{equation*}
Here $J_{l}^{p+1}(\xi,\eta)$ is the Jacobian of the mapping
$M_{l}^{p+1}$ from $S$ to $\Omega_{l}^{p+1}$. Let
$\mathcal{L}_{l}^{p+1}(\xi,\eta)=\mathcal{L}(\xi,\eta)\,\sqrt{J_{l}^{p+1}}.$
Once more we can define ($\mathcal{L}_{l}^{p+1})^{a},$ a
differential operator which is an approximation to
$\mathcal{L}_{l}^{p+1}$ in which the coefficients of
$\mathcal{L}_{l}^{p+1}$ are replaced by polynomial approximations.
It can be shown as before that
\begin{align*}
&\sum_{l=1}^{L}\int_{\Omega_{l}^{p+1}}|\mathcal{L}w|^{2}\,
\hbox{d}x_{1}\hbox{d}x_{2}\\[.4pc]
&\quad\, \leq
C\,\sum_{l=1}^{L}\int_{S}|(\mathcal{L}_{l}^{p+1})^{^{a}}w_{l}^{p+1}(\xi,\eta)|^{2}\,
\hbox{d}\xi \hbox{d}\eta+\varepsilon_{_{W}}\sum_{l=1}^{L}\Vert
w_{l}^{p+1}(\xi,\eta)\Vert_{_{2,S}}^{2}.
\end{align*}
Here $C$ is a constant and $\varepsilon_{_{W}}\rightarrow0$ as
$W\rightarrow$$\infty.$ In fact, $\varepsilon_{_{W}}$ is
exponentially small in $W$.

We now prove a result which we shall need in the
sequel.

\begin{lem}
Let $\omega\in H_{\beta}^{2,2}(\Omega)$. Then there exists a
constant $C$ such that
\begin{align}
&\frac{1}{C}\Bigg(\sum_{k=1}^{p}\Bigg(|\omega(A_{k})|^{2}+\sum_{|\alpha|\leq2}\int_{\hat{\Omega}^{k}}|D_{\nu_{k},\phi_{k}}^{\alpha}(\omega(\nu_{k},\phi_{k})-\omega(A_{k}))|^{2} \nonumber\\[.4pc]
&\qquad\, \times {\rm e}^{-2\lambda_{k}\nu_{k}}{\rm d}\nu_{k}{\rm
d}\phi_{k}\Bigg)+\Vert
\omega(x_{1},x_{2})\Vert_{_{H^{2}(\Omega^{^{p+1}})}}^{2}\Bigg)\nonumber
\end{align}
\begin{align}
&\quad\, \leq \Vert \omega\Vert_{_{H_{\beta}^{2,2}(\Omega)}}^{2}\nonumber\\[.4pc]
&\quad\, \leq C\Bigg(\sum_{k=1}^{p}\Bigg(|\omega(A_{k})|^{2}+
\sum_{|\alpha|\leq2}\int_{\hat{\Omega}^{k}}
|D_{\nu_{k},\phi_{k}}^{\alpha}(\omega(\nu_{k},\phi_{k})-\omega(A_{k}))|^{2}\nonumber\\[.4pc]
&\qquad\, \times {\rm e}^{-2\lambda_{k}\nu_{k}}{\rm d}\nu_{k}{\rm
d}\phi_{k}\Bigg)+\Vert \omega(x_{1},x_{2})\Vert_{_{H^{2}
(\Omega^{^{p+1}})}}^{2}\Bigg).\label{1..9}
\end{align}
\end{lem}

Here $\lambda_{k}=1-\beta_{k}$.

\begin{proof}
Let $\psi_{_{k}}\in C_{0}^{^{\infty}}(R)$ such that
$\psi_{_{k}}(r_{k})=1$ for $r_{k}\leq\rho$ and
$\psi_{_{k}}(r_{k})=0$ for $r_{k}\geq\rho^{1}$ for
$k=1,2,\dots,p.$ Here $\rho^{1}>$~$\rho$ is chosen so that
$\Omega_{\rho^{1}}^{k}=\{(x_{1},x_{2})\hbox{\rm :}\
r_{k}\leq\rho^{1}\}$ have the property that
$\Omega_{\rho^{1}}^{k}\cap\Omega_{\rho^{1}}^{l}=\emptyset$ if
$k\neq l.$ We define $\omega_{k}=\omega\,\psi_{k}$ for
$k=1,\dots,p$ and $\omega_{0}=1-{\sum_{k=1}^{p}\omega_{k}}.$ Then
$\omega_{k}\in H_{_{\beta}}^{^{2,2}}(\Omega)$ for $k=1,\dots,p.$

Now by Lemma~2.1 of \cite{babguo2},
$H_{_{\beta}}^{^{2,2}}(\Omega)\subseteq C(\bar{\Omega})$ with
continuous injection. Hence we conclude that
\begin{equation*}
\sum_{k=1}^{p}|\omega_{k}(A_{k})|^{2}\leq C\,\,\sum_{k=1}^{p}\Vert
\omega_{k}\Vert _{_{H_{\beta}^{2,2}(\Omega)}}^{2}.
\end{equation*}
Therefore
\begin{equation}
\sum_{k=1}^{p}|\omega(A_{k})|^{2}\leq C\,\,\sum_{k=1}^{p}\Vert
\omega\Vert _{_{H_{\beta}^{2,2}(\Omega)}}^{2}.\label{1.9}
\end{equation}
We now cite Lemma~2.2 of \cite{babguo3}. Let $u\in
H_{_{\beta}}^{^{2,2}}(\Omega)$. Then
\begin{enumerate}
\renewcommand{\labelenumi}{(\roman{enumi})}
\leftskip .2pc
\item
\begin{equation*}
\hskip -1.25pc \sum_{|\alpha|=1}\Vert
D^{\alpha}u\,\Phi_{\beta-1}\Vert_{_{L^{2}(\Omega)}}\leq C\Vert
u\Vert_{_{H_{\beta}^{2,2}(\Omega)}}.
\end{equation*}

\item Let $u(A_{i})=0$, for $i=1,\dots,p$. Then
\begin{equation*}
\hskip -1.25pc \Vert u\,\Phi_{\beta-2}\Vert_{_{L^{2}(\Omega)}}\leq
C\Vert u\Vert_{_{H_{\beta}^{2,2}(\Omega)}}.
\end{equation*}
\end{enumerate}\vspace{-.6pc}
From (i) we obtain
\begin{equation}
\sum_{k=1}^{p}\int_{\hat{\Omega}^{k}}\sum_{|\alpha|=1}|D_{\nu_{k},\phi_{k}}^{\alpha}\omega(\nu_{k},\phi_{k})|^{2}\hbox{e}^{-2(1-\beta_{k})\nu_{k}}\,
\hbox{d}\nu_{k}\hbox{d}\phi_{k}\leq C(\left\Vert \omega\right\Vert
_{_{H_{\beta}^{2,2}(\Omega)}}^{2}).\label{1.10}
\end{equation}
Here $C$ is a generic constant. Now using (ii) we get
\begin{align*}
&\int_{\Omega}|\omega_{k}(\nu_{k},\phi_{k})-\omega(A_{k})\psi_{k}|^{2}\hbox{e}^{-2(1-\beta_{k})\nu_{k}}\,
\hbox{d}\nu_{k}\hbox{d}\phi_{k}\\[.4pc]
&\quad\, \leq C(\Vert \omega_{k}\Vert
_{_{H_{\beta}^{2,2}(\Omega)}}^{2}+|\omega(A_{k})|^{2}).
\end{align*}
Hence
\begin{equation}
\sum_{k=1}^{p}\int_{\hat{\Omega}^{k}}|\omega(\nu_{k},\phi_{k})-\omega(A_{k})|^{2}\hbox{e}^{-2(1-\beta_{k})\nu_{k}}\,
\hbox{d}\nu_{k}\hbox{d}\phi_{k}\leq C\Vert
\omega\Vert_{_{H_{\beta}^{2,2}(\Omega)}}^{2}.\label{1.11}
\end{equation}
Finally,
\begin{equation}
\sum_{k=1}^{p}\int_{\hat{\Omega}^{k}}\sum_{|\alpha|=2}|D_{\nu_{k},\phi_{k}}^{\alpha}\omega(\nu_{k},\phi_{k})|^{2}\hbox{e}^{-2(1-\beta_{k})\nu_{k}}\,
\hbox{d}\nu_{k}\hbox{d}\phi_{k}\leq C \Vert
\omega\Vert_{_{H_{\beta}^{2,2}(\Omega)}}^{2}.\label{1.12}
\end{equation}
Combining the estimates (3.10)--(3.13) we get (3.9).\hfill
$\blacksquare$
\end{proof}

We now introduce some notation which is needed to state the
stability estimate Theorem 3.1 which is the main result of this
section.

Let $\gamma_{s}$ be a side common to the elements
$\Omega_{m}^{p+1}$ and $\Omega_{n}^{p+1}$ and let
$\gamma_{s}\subseteq\Omega^{p+1}.$ We may assume that $\gamma_{s}$
is the image of $\eta=-1$ under the mapping $M_{m}^{p+1}$ which
maps $S$ to $\Omega_{m}^{p+1}$and also the image of $\eta=1$ under
the mapping $M_{n}^{p+1}$ which maps $S$ to $\Omega_{n}^{p+1}.$ By
the chain rule
\begin{align*}
(u_{m}^{p+1})_{x_{1}}
&=(u_{m}^{p+1})_{\xi}\,\,\xi_{x_{1}}+(u_{m}^{p+1})_{\eta}\,\,\eta_{x_{1}},\quad\textrm{and}\\[.4pc]
(u_{m}^{p+1})_{x_{2}}
&=(u_{m}^{p+1})_{\xi}\,\,\xi_{x_{2}}+(u_{m}^{p+1})_{\eta}\,\,\eta_{x_{2}}.
\end{align*}
Now let $\hat{\xi}_{x_{1}}$ denote the polynomial approximation of
$\xi_{x_{1}}(\xi,\eta),$ of degree $W$ in $\xi$ and $\eta$
separately, as defined in Theorem 4.46 of \cite{schwab}. In the
same way $\hat{\eta}_{x_{1}}$,$\,\hat{\xi}_{x_{2}}$ and
$\hat{\eta}_{x_{2}}$ can be defined. We now
define
\begin{align*}
(u_{m}^{p+1})_{x_{1}}^{a}&=(u_{m}^{p+1})_{\xi}\,\,\hat{\xi}_{x_{1}}+(u_{m}^{p+1})_{\eta}\,\,\hat{\eta}_{x_{1}},\quad\textrm{and}\\[.4pc]
(u_{m}^{p+1})_{x_{2}}^{a}&=(u_{m}^{p+1})_{\xi}\,\,\hat{\xi}_{x_{2}}+(u_{m}^{p+1})_{\eta}\,\,\hat{\eta}_{x_{2}}.
\end{align*}
Let
\begin{align*}
\Vert [u^{p+1}]\Vert_{_{0,\gamma_{s}}}^{2}&=\Vert
u_{m}^{p+1}(\xi,-1)-u_{n}^{p+1}(\xi,1)\Vert_{_{0,I}}^{2},\\[.4pc]
\Vert [(u_{x_{1}}^{p+1})^{a}]\Vert_{_{1/2,\gamma_{s}}}^{2}&=\Vert (u_{m}^{p+1})_{x_{1}}^{a}(\xi,-1)-(u_{n}^{p+1})_{x_{1}}^{a}(\xi,1)\Vert_{_{1/2,I}}^{2},\quad\textrm{and}\\[.4pc]
\Vert [(u_{x_{2}}^{p+1})^{a}]\Vert _{_{1/2,\gamma_{s}}}^{2}&=\Vert
(u_{m}^{p+1})_{x_{2}}^{a}(\xi,-1)-(u_{n}^{p+1})_{x_{2}}^{a}(\xi,1)\Vert_{_{1/2,I}}^{2}.
\end{align*}
Here $I=(-1,1).$ Next, let
$\gamma_{s}\subseteq\Gamma^{[0]}\cap\partial\Omega^{p+1}$ and let
$\gamma_{s}$ be the image of $\eta=-1$ under the mapping
$M_{m}^{p+1}$ which maps $S$ to $\Omega_{m}^{p+1}.$ We can define
$\big(\frac{\partial u_{m}^{p+1}}{\partial T}\big)^{a},$ an
approximation to $\frac{\partial u^{p+1}}{\partial T}$ as before.
Let
\begin{align*}
&\Vert u^{p+1}\Vert_{_{0,\gamma_{s}}}^{2}+\left\Vert
\left(\frac{\partial u^{p+1}}{\partial T}\right)^{a}\right\Vert
_{_{1/2,\gamma_{s}}}^{2}\\[.4pc]
&\quad\, =\Vert u_{m}^{p+1}(\xi,-1)\Vert _{_{0,I}}^{2}+\left\Vert
\left(\frac{\partial u_{m}^{p+1}}{\partial
T}\right)^{a}(\xi,-1)\right\Vert _{_{1/2,I}}^{2}.
\end{align*}
In the same way, if
$\gamma_{s}\subseteq\Gamma^{[1]}\cap\partial\Omega^{p+1},$
$\big\Vert \big(\frac{\partial u^{p+1}}{\partial
N}\big)_{A}^{a}\big\Vert_{_{1/2,\gamma_{s}}}^{2}$ can be
defined.\pagebreak

Let $\gamma_{s}\subseteq\Gamma^{[1]}\cap\partial\Omega^{k}$ for
$1\leq k\leq p.$ Let $\tilde{\gamma}_{s}$ be the image of
$\gamma_{s}$ in $(\tau_{k},\theta_{k})$ coordinates and
$\hat{\gamma}_{s}$ be the image of $\gamma_{s}$ in
$(\nu_{k},\phi_{k})$ coordinates. Let $(n_{1},n_{2})$ be the
normal at a point $\tilde{p}$ on $\tilde{\gamma}_{s}.$ Define
\begin{equation*}
\left(\frac{\partial u^{k}}{\partial
n}\right)_{\tilde{A}^{k}}=\sum_{i,j=1}^{2}n_{i}\,\tilde{a}_{i,j}^{k}\,\frac{\partial
u^{k}}{\partial y_{j}}.
\end{equation*}
Now $\hat{\gamma}_{s}$ is a portion of the straight line
$\phi_{k}=\alpha$, where $\alpha$ is a constant. Let
$\big(\frac{\partial u^{k}}{\partial n}\big)_{\tilde{A}^{k}}^{a}$,
denote an approximation to $\big(\frac{\partial u^{k}}{\partial
n}\big)_{\tilde{A}^{k}}$ as before, and using this $\big\Vert
\big(\frac{\partial u^{k}}{\partial
n}\big)_{\tilde{A}^{k}}^{a}\big\Vert_{_{1/2,\hat{\gamma}_{s}}}^{2}$
can be defined. Let $\gamma_{s}\subseteq\bar{\Omega}^{k}.$ Define
\begin{equation*}
d(A_{k},\gamma_{s})=\inf_{_{x\in\gamma_{s}}}\{
\hbox{distance}(A_{k},x)\}. \end{equation*} Let
\begin{align}
&\mathcal{V}_{_{\rm vertices}}^{^{M,W}}(\{ u_{i,j}^{k}(\nu_{k},\phi_{k})\}_{i,j,k},\{
u_{l}^{_{^{p+1}}}(\xi,\eta)\}_{l})\nonumber\\[.2pc]
&\quad\, =
\sum_{k=1}^{p}\sum_{j=2}^{M}\sum_{i=1}^{I_{k}}(\rho\mu_{k}^{M+1-j})^{-2\lambda_{k}}\Vert
 (\mathcal{L}_{i,j}^{k})^{^{a}}u_{i,j}^{k}(\nu_{k},\phi_{k})\Vert_{_{^{0,\hat{\Omega}_{i,j}^{k}}}}^{2}\nonumber\\
&\qquad\, + \sum_{k=1}^{p}\sum_{\gamma_{s}\subseteq\Omega^{k}\cup
B_{\rho}^{^{k}},\mu(\hat{\gamma}_{s})<\infty}
d(A_{k},\gamma_{s})^{-2\lambda_{k}}\nonumber\\[.4pc]
&\qquad\, \times (\Vert
[u^{k}]\Vert_{_{0,\hat{\gamma}_{s}}}^{2}+\Vert
[(u_{\nu_{k}}^{k})^{a}]\Vert_{_{1/2,\hat{\gamma}_{s}}}^{2}+\Vert
 [(u_{\phi_{k}}^{k})^{a}]\Vert_{_{1/2,\hat{\gamma}_{s}}}^{2})\nonumber\\
&\qquad\, +
\sum_{l\in\mathcal{D}}\sum_{k=l-1}^{l}(|h_{k}|^{2}+\sum_{\gamma_{s}\subseteq\partial\Omega^{k}\cap\Gamma_{l},\mu(\hat{\gamma}_{s})<\infty}
d(A_{k},\gamma_{s})^{-2\lambda_{k}}\nonumber\\[.4pc]
&\qquad\, \times (\Vert
u^{k}-h_{k}\Vert_{_{0,\hat{\gamma}_{s}}}^{2}+\Vert
 u_{\nu_{k}}^{k}\Vert_{_{1/2,\hat{\gamma}_{s}}}^{2}))\nonumber\\
&\qquad\, +
\sum_{l\in\mathcal{N}}\sum_{k=l-1}^{l}\sum_{\gamma_{s}\subseteq\partial\Omega^{k}\cap\Gamma_{l},\mu(\hat{\gamma}_{s})<\infty}
d(A_{k},\gamma_{s})^{-2\lambda_{k}}\left\Vert \left(\frac{\partial
u^{k}}{\partial
n}\right)_{\tilde{A}^{k}}^{a}\right\Vert_{_{1/2,\hat{\gamma}_{s}}}^{2}.\label{1.14a}
\end{align}
Here $\{ \{ u_{i,j}^{k}(\nu_{k},\phi_{k})\} _{i,j,k},\{
u_{l}^{_{^{p+1}}}(\xi,\eta)\}_{l}\} \in\Pi^{M,W}$ and
$u_{i,1}^{k}=h_{k}$ for $1\leq i\leq I_{k}.$ Moreover
$\mu(\hat{\gamma}_{s})$ denotes the measure of $\hat{\gamma}_{s}.$
Next, we define
\begin{align*}
& \mathcal{V}_{_{\rm interior}}^{^{M,W}}(\{ u_{i,j}^{k}(\nu_{k},\phi_{k})\}_{i,j,k},\{
u_{l}^{_{^{p+1}}}(\xi,\eta)\}_{l})\\[.4pc]
&\quad\, = \sum_{l=1}^{L}\Vert (\mathcal{L}_{l}^{p+1})^{^{a}}u_{l}^{p+1}(\xi,\eta)\Vert_{_{0,S}}^{2}\\[.4pc]
 &\qquad\, +\sum_{\gamma_{s}\subseteq\Omega^{^{p+1}}}(\Vert [u^{p+1}]\Vert_{_{0,\gamma_{s}}}^{2}+\Vert [(u_{x_{1}}^{p+1})^{a}]\Vert_{_{1/2,\gamma_{s}}}^{2}+\Vert [(u_{x_{2}}^{p+1})^{a}]\Vert_{_{1/2,\gamma_{s}}}^{2})
 \end{align*}
 \begin{align*}
 &\qquad\, + \sum_{l\in\mathcal{D}}\sum_{\gamma_{s}\subseteq\partial\Omega^{^{p+1}}\cap\Gamma_{l}}\left(\Vert u^{p+1}\Vert_{_{0,\gamma_{s}}}^{2}+\left\Vert \left(\frac{\partial u^{p+1}}{\partial T}\right)^{a}\right\Vert_{_{1/2,\gamma_{s}}}^{2}\right)\\[.4pc]
 &\qquad\, + \sum_{l\in\mathcal{N}}\sum_{\gamma_{s}\subseteq\partial\Omega^{^{p+1}}\cap\Gamma_{l}}\left\Vert \left(\frac{\partial u^{p+1}}{\partial N}\right)_{A}^{a}\right\Vert_{_{1/2,\gamma_{s}}}^{2}.
 \end{align*}
Let
\begin{align}
 &\mathcal{V}^{^{M,W}}(\{ u_{i,j}^{k}(\nu_{k},\phi_{k})\}_{i,j,k},\{
 u_{l}^{_{^{p+1}}}(\xi,\eta)\}_{l})\nonumber\\[.4pc]
 &\quad\,=  \mathcal{V}_{_{\rm vertices}}^{^{M,W}}(\{ u_{i,j}^{k}(\nu_{k},\phi_{k})\}_{i,j,k},\{ u_{l}^{_{^{p+1}}}(\xi,\eta)\}_{l})\nonumber\\[.4pc]
 &\qquad\, +\mathcal{V}_{_{\rm interior}}^{^{M,W}}(\{ u_{i,j}^{k}(\nu_{k},\phi_{k})\}_{i,j,k},\{ u_{l}^{_{^{p+1}}}(\xi,\eta)\}_{l}).\label{1.15}
 \end{align}
We can now state the main result of this section.

\begin{theo}[\!]
For $M$ and $W$ large enough the estimate
\begin{align}
& \sum_{k=1}^{p}\left(\left|h_{k}\right|^{2}+\sum_{i=1}^{I_{k}}\sum_{j=2}^{M}(\rho\mu_{k}^{M+1-j})^{-2\lambda_{k}}\Vert
u_{i,j}^{k}(\nu_{k},\phi_{k})-h_{k}\Vert_{_{2,\hat{\Omega}_{i,j}^{k}}}^{2}\right)\nonumber\\[.4pc]
&\qquad\, +\sum_{l=1}^{L}\Vert
u_{l}^{p+1}(\xi,\eta)\Vert_{_{2,S}}^{2}\nonumber\\[.4pc]
&\quad\, \leq C(\ln W)^{2}\,\mathcal{V}^{^{M,W}}(\{
u_{i,j}^{k}(\nu_{k},\phi_{k})\}_{i,j,k},\{
u_{l}^{_{^{p+1}}}(\xi,\eta)\}_{l})\label{1.16}
\end{align}
holds. Here $C$ is a constant.
\end{theo}

\begin{proof}
By Lemma 7.1 there exist $\{ \{
v_{i,j}^{k}(\nu_{k},\phi_{k})\}_{i,j,k},\{
v_{l}^{_{^{p+1}}}(\xi,\eta)\}_{l}\} $ such that $w$ defined as
$w=u+v\in H_{_{\beta}}^{^{2,2}}(\Omega).$ Moreover $v_{i,1}^{k}=0$
for all $i$ and $k.$ Hence by Theorem 2.1 of \cite{babguo1},
\begin{equation}
\left\Vert w\right\Vert
_{_{H_{_{\beta}}^{^{2,2}}(\Omega)}}^{2}\leq C\left(\left\Vert
\mathcal{L}w\right\Vert_{_{L_{\beta}(\Omega)}}^{2}+\left\Vert
w\right\Vert
_{_{H_{_{\beta}}^{^{\frac{3}{2},\frac{3}{2}}}(\Gamma^{[0]})}}^{2}+\left\Vert
\left(\frac{\partial w}{\partial N}\right)_{A}\right\Vert
_{_{H_{_{\beta}}^{^{\frac{1}{2},\frac{1}{2}}}(\Gamma^{[1]})}}^{2}\right).\label{1.17}
\end{equation}
Now $v_{i,1}^{k}(\nu_{k},\phi_{k})=0$ for $1\leq i\leq I_{k}$.
Hence by (3.8),
\begin{align*}
\hskip -4pc &\left\Vert \mathcal{L}w\right\Vert_{_{L_{_{\beta}}(\Omega)}}^{2}\nonumber\\[.4pc]
\hskip -4pc &\quad\, \leq
2\left(\sum_{k=1}^{p}\sum_{j=2}^{M}\sum_{i=1}^{I_{k}}\Vert
(\mathcal{L}_{i,j}^{k})^{^{a}}u_{i,j}^{k}(\nu_{k},\phi_{k})\Vert_{_{0,\hat{\Omega}_{i,j}^{k}}}^{2}
+\sum_{l=1}^{L}\Vert
(\mathcal{L}_{l}^{p+1})^{^{a}}u_{l}^{p+1}(\xi,\eta)\Vert_{_{0,S}}^{2}\right)\nonumber
\end{align*}
\begin{align}
&\hskip -2.5pc +  C\left(\sum_{k=1}^{p}\sum_{i=1}^{I_{k}}\sum_{j=2}^{M}(\rho\mu_{k}^{M+1-j})^{-2\lambda_{k}}\Vert v_{i,j}^{k}(\nu_{k},\phi_{k})\Vert_{_{2,\hat{\Omega}_{i,j}^{k}}}^{2} +\sum_{l=1}^{L}\Vert v_{l}^{p+1}(\xi,\eta)\Vert_{_{2,S}}^{2}\right)\nonumber\\[.4pc]
&\hskip -2.5pc  + \varepsilon_{_{W}}\left(\sum_{k=1}^{p}\sum_{j=2}^{M}\sum_{i=1}^{I_{k}}(\rho\mu_{k}^{M+1-j})^{-2\lambda_{k}}\Vert u_{i,j}^{k}(\nu_{k},\phi_{k})-h_{k}\Vert_{_{2,\hat{\Omega}_{i,j}^{k}}}^{2} +\sum_{k=1}^{p}\left|h_{k}\right|^{2}\right)\nonumber\\[.4pc]
&\hskip -2.5pc  + \varepsilon_{_{W}}\left(\sum_{l=1}^{L}\Vert
u_{l}^{p+1}(\xi,\eta)\Vert_{_{2,S}}^{2}\right)+\varepsilon_{_{M}}\left(\sum_{k=1}^{p}\left|h_{k}\right|^{2}\right).\label{3.18}
\end{align}
Now using Lemma 7.2,
\begin{align}
 &\Vert w\Vert_{H_{_{\beta}}^{^{\frac{3}{2},\frac{3}{2}}}(\Gamma^{[0]})}^{2}+\left\Vert \left(\frac{\partial w}{\partial
 N}\right)_{A}\right\Vert_{H_{_{\beta}}^{^{\frac{1}{2},\frac{1}{2}}}(\Gamma^{[1]})}^{2}\nonumber\\[.4pc]
 &\quad\, \leq  C\,(\ln\, W)^{2}\left(\sum_{_{k\hbox{\rm :}\ \partial\Omega^{k}\cap\Gamma^{[0]}\neq\emptyset}}^{p}|h_{k}|^{2}+\sum_{l\in\mathcal{D}}\sum_{k=l-1}^{l}\right.\nonumber\\[.4pc]
 &\qquad\, \times \sum_{\gamma_{s}\subseteq\partial\Omega^{k}\cap\Gamma_{l},\mu(\hat{\gamma}_{s})<\infty} d(A_{k},\gamma_{s})^{-2\lambda_{k}}(\Vert u^{k}-h_{k}\Vert_{0,\hat{\gamma}_{s}}^{2}+\Vert u_{\nu_{k}}^{k}\Vert_{1/2,\hat{\gamma}_{s}}^{2})\nonumber\\[.4pc]
 &\qquad\, +  \sum_{l\in\mathcal{N}}\sum_{k=l-1}^{l}\sum_{\gamma_{s}\subseteq\partial\Omega^{k}\cap\Gamma_{l},\mu(\hat{\gamma}_{s})<\infty} d(A_{k},\gamma_{s})^{-2\lambda_{k}}\left\Vert \left(\frac{\partial u^{k}}{\partial n}\right)_{\tilde{A}^{k}}^{a}\right\Vert_{1/2,\hat{\gamma}_{s}}^{2}\nonumber \\[.4pc]
 &\qquad\, +  \sum_{l\in\mathcal{D}}\sum_{\gamma_{s}\subseteq\partial\Omega^{^{p+1}}\cap\Gamma_{l}}\left(\Vert u^{p+1}\Vert_{_{0,\gamma_{s}}}^{2}+\left\Vert \left(\frac{\partial u^{p+1}}{\partial T}\right)^{a}\right\Vert_{_{1/2,\gamma_{s}}}^{2}\right)\nonumber\\[.4pc]
 &\qquad\, +  \sum_{l\in\mathcal{N}}\sum_{\gamma_{s}\subseteq\partial\Omega^{^{p+1}}\cap\Gamma_{l}}\left\Vert \left(\frac{\partial u^{p+1}}{\partial N}\right)_{A}^{a}\right\Vert_{_{1/2,\gamma_{s}}}^{2}\nonumber \\[.4pc]
 &\qquad\, +  \sum_{k=1}^{p}\sum_{\gamma_{s}\subseteq\Omega^{k}\cup B_{\rho}^{k},\mu(\hat{\gamma}_{s})<\infty} d(A_{k},\gamma_{s})^{-2\lambda_{k}}\nonumber\\[.4pc]
 &\qquad\, \times (\Vert [u^{k}]\Vert_{_{0,\hat{\gamma}_{s}}}^{2}+\Vert [(u_{\nu_{k}}^{k})^{a}]\Vert_{_{1/2,\hat{\gamma}_{s}}}^{2}+\Vert [(u_{\phi_{k}}^{k})^{a}]\Vert_{_{1/2,\hat{\gamma}_{s}}}^{2})\nonumber \\[.4pc]
 &\qquad\, +  \sum_{\gamma_{s}\subseteq\Omega^{^{p+1}}}(\Vert [u^{p+1}]\Vert_{_{0,\gamma_{s}}}^{2}+\Vert [(u_{x_{1}}^{p+1})^{a}]\Vert_{_{1/2,\gamma_{s}}}^{2}+\Vert [(u_{x_{2}}^{p+1})^{a}]\Vert_{_{1/2,\gamma_{s}}}^{2}))\nonumber \\[.4pc]
 &\qquad\, +  \varepsilon_{_{W}}\,\left(\sum_{k=1}^{p}|h_{k}|^{2}+\sum_{k=1}^{p}\sum_{i=1}^{I_{k}}\sum_{j=2}^{M}(\rho\mu_{k}^{M+1-j})^{-2\lambda_{k}}\Vert u_{i,j}^{k}(\nu_{k},\phi_{k}) \right.\nonumber\\[.4pc]
 &\qquad\, \left. -h_{k}\Vert_{_{2,\hat{\Omega}_{i,j}^{k}}}^{2}+\sum_{l=1}^{L}\Vert u_{l}^{p+1}(\xi,\eta)\Vert_{_{2,S}}^{2}\right).\label{3.19}
 \end{align}
Combining (3.17)--(3.19) we obtain
\begin{align}
\hskip -4pc \Vert w\Vert_{_{H_{_{\beta}}^{^{2,2}}(\Omega)}}^{2} & \leq C\,(\ln W)^{2}\,\mathcal{V}^{^{M,W}}(\{ u_{i,j}^{k}(\nu_{k},\phi_{k})\}_{i,j,k},\{
u_{l}^{_{^{p+1}}}(\xi,\eta)\}_{l})\nonumber\\[.4pc]
\hskip -4pc  &\ \ \ \, +\!C\!\left(\sum_{k=1}^{p}\sum_{i=1}^{I_{k}}\sum_{j=2}^{M}(\rho\mu_{k}^{M+1-j})^{-2\lambda_{k}}\Vert v_{i,j}^{k}(\nu_{k},\phi_{k})\Vert_{_{2,\hat{\Omega}_{i,j}^{k}}}^{2}\!\!+\!\!\sum_{l=1}^{L}\Vert v_{l}^{p+1}(\xi,\eta)\Vert_{_{2,S}}^{2}\!\right)\nonumber\\[.4pc]
\hskip -4pc  &\ \ \ \,+\varepsilon_{_{W}}\left(\sum_{k=1}^{p}\sum_{j=2}^{M}\sum_{i=1}^{I_{k}}(\rho\mu_{k}^{M+1-j})^{-2\lambda_{k}}\Vert u_{i,j}^{k}(\nu_{k},\phi_{k})-h_{k}\Vert_{_{2,\hat{\Omega}_{i,j}^{k}}}^{2}+\sum_{k=1}^{p}|h_{k}|^{2}\right)\nonumber \\[.4pc]
\hskip -4pc  &\ \ \ \,
+\varepsilon_{_{W}}\left(\sum_{l=1}^{L}\Vert
u_{l}^{p+1}(\xi,\eta)\Vert
_{_{2,S}}^{2}\right)+\varepsilon_{_{M}}\left(\sum_{k=1}^{p}|h_{k}|^{2}\right).\label{3.20}
\end{align}
Now using (3.9),
\begin{align}
&\sum_{k=1}^{p}\left(\left|h_{k}\right|^{2}+\sum_{i=1}^{I_{k}}\sum_{j=2}^{M}(\rho\mu_{k}^{M+1-j})^{-2\lambda_{k}}\Vert u_{i,j}^{k}(\nu_{k},\phi_{k})-h_{k}\Vert_{_{2,\hat{\Omega}_{i,j}^{k}}}^{2}\right)\nonumber\\[.4pc]
&\qquad\, +\sum_{l=1}^{L}\Vert u_{l}^{p+1}(\xi,\eta)\Vert_{_{2,S}}^{2}\nonumber\\[.4pc]
&\quad\, \leq  K\left(\Vert
w\Vert_{_{H_{_{\beta}}^{^{2,2}}(\Omega)}}^{2}+\sum_{k=1}^{p}\sum_{i=1}^{I_{k}}\sum_{j=2}^{M}(\rho\mu_{k}^{M+1-j})^{-2\lambda_{k}}\Vert
v_{i,j}^{k}(\nu_{k},\phi_{k})\Vert_{_{^{2,\hat{\Omega}_{i,j}^{k}}}}^{2}\right.\nonumber\\[.4pc]
&\qquad\, \left. +\sum_{l=1}^{L}\Vert v_{l}^{p+1}(\xi,\eta)\Vert
_{_{2,S}}^{2}\right).\label{3.21}
 \end{align}
Combining (3.20) and (3.21) gives
\begin{align}
&\hskip -4pc
\sum_{k=1}^{p}\left(|h_{k}|^{2}+\sum_{i=1}^{I_{k}}\sum_{j=2}^{M}(\rho\mu_{k}^{M+1-j})^{-2\lambda_{k}}\Vert
u_{i,j}^{k}(\nu_{k},\phi_{k})-h_{k}\Vert_{_{2,\hat{\Omega}_{i,j}^{k}}}^{2}\right)+\sum_{l=1}^{L}\Vert
u_{l}^{p+1}(\xi,\eta)\Vert_{_{2,S}}^{2}\nonumber\\[.4pc]
 &\hskip -3pc  \leq  C\,(\ln W)^{2}\,\mathcal{V}^{^{M,W}}(\{ u_{i,j}^{k}(\nu_{k},\phi_{k})\}_{i,j,k},\{ u_{l}^{_{^{p+1}}}(\xi,\eta)\}_{l})\nonumber \\[.4pc]
&\hskip -2pc  + C\left(\sum_{k=1}^{p}\sum_{i=1}^{I_{k}}\sum_{j=2}^{M}(\rho\mu_{k}^{M+1-j})^{-2\lambda_{k}}\Vert v_{i,j}^{k}(\nu_{k},\phi_{k})\Vert_{_{2,\hat{\Omega}_{i,j}^{k}}}^{2}+\sum_{l=1}^{L}\Vert v_{l}^{p+1}(\xi,\eta)\Vert_{_{2,S}}^{2}\right)\nonumber \\[.4pc]
&\hskip -2pc   +  \varepsilon_{_{W}}\left(\sum_{k=1}^{p}\sum_{j=2}^{M}\sum_{i=1}^{I_{k}}(\rho\mu_{k}^{M+1-j})^{-2\lambda_{k}}\Vert u_{i,j}^{k}(\nu_{k},\phi_{k})-h_{k}\Vert_{_{2,\hat{\Omega}_{i,j}^{k}}}^{2}+\sum_{k=1}^{p}|h_{k}|^{2}\right)\nonumber\\[.4pc]
 &\hskip -2pc+\varepsilon_{_{W}}\left(\sum_{l=1}^{L}\Vert
u_{l}^{p+1}(\xi,\eta)\Vert_{_{2,S}}^{2}\right)+\varepsilon_{_{M}}\left(\sum_{k=1}^{p}\left|h_{k}\right|^{2}\right).\label{3.22}
\end{align}
Now by Lemma 7.1,
\begin{align}
&\hskip -4pc
\sum_{k=1}^{p}\sum_{i=1}^{I_{k}}\sum_{j=2}^{M}(\rho\mu_{k}^{M+1-j})^{-2\lambda_{k}}\Vert
v_{i,j}^{k}(\nu_{k},\phi_{k})\Vert_{_{2,\hat{\Omega}_{i,j}^{k}}}^{2}+\sum_{l=1}^{L}\Vert
v_{l}^{p+1}(\xi,\eta)\Vert_{_{2,S}}^{2}\nonumber\\[.4pc]
&\hskip -3pc  \leq  C\,(\ln W)^{2}\left(\sum_{k=1}^{p}\sum_{\gamma_{s}\subseteq\Omega^{k}\cup B_{\rho}^{k},\mu(\hat{\gamma}_{s})<\infty} d(A_{k},\gamma_{s})^{-2\lambda_{k}}\right.\nonumber\\[.4pc]
&\hskip -2pc  \times (\Vert [u^{k}]\Vert_{_{0,\hat{\gamma}_{s}}}^{2}+\Vert [(u_{\nu_{k}}^{k})^{a}]\Vert_{_{1/2,\hat{\gamma}_{s}}}^{2}+\Vert [(u_{\phi_{k}}^{k})^{a}]\Vert_{_{1/2,\hat{\gamma}_{s}}}^{2})\nonumber\\[.4pc]
&\hskip -2pc  + \left.\sum_{\gamma_{s}\subseteq\Omega^{^{p+1}}}(\Vert [u^{p+1}]\Vert_{_{0,\gamma_{s}}}^{2}+\Vert [(u_{x_{1}}^{p+1})^{a}]\Vert_{_{1/2,\gamma_{s}}}^{2}+\Vert [(u_{x_{2}}^{p+1})^{a}]\Vert_{_{1/2,\gamma_{s}}}^{2})\right)\nonumber \\
&\hskip -2pc  +
\varepsilon_{_{W}}\left(\sum_{k=1}^{p}\sum_{i=1}^{I_{k}}\sum_{j=2}^{M}(\rho\mu_{k}^{M+1-j})^{-2\lambda_{k}}\Vert
u_{i,j}^{k}(\nu_{k},\phi_{k})-h_{k}\Vert
_{_{2,\hat{\Omega}_{i,j}^{k}}}^{2}\right.\nonumber\\[.4pc]
&\hskip -2pc \left.+\sum_{l=1}^{L}\Vert
u_{l}^{p+1}(\xi,\eta)\Vert_{_{2,S}}^{2}\right).\label{3.23}
\end{align}
Combining (3.22) and (3.23) we get the result.\hfill
$\blacksquare$
\end{proof}

\section{The numerical scheme}

\setcounter{equation}{0}

As in \S3,
\begin{equation*}
\hat{\Omega}_{i,j}^{k}=\{ (\nu_{k},\phi_{k})\hbox{\rm :}\
\nu_{j}^{k}<\nu_{k}<\nu_{j+1}^{k},\psi_{i}^{k}<\phi_{k}<\psi_{i+1}^{k}\}
\end{equation*}
for $1\leq j\leq M,1\leq i\leq I_{k,j},1\leq k\leq p$ in $\nu_{k}$
and $\phi_{k}$ variables.

We now define a nonconforming spectral element representation on
each of these subdomains as follows:
\begin{equation*}
u_{i,j}^{k}(\nu_{k},\phi_{k})=h_{k},\quad \mathrm{if}\: j=1,1\leq
i\leq I_{k},1\leq k\leq p
\end{equation*}
and
\begin{equation*}
u_{i,j}^{k}(\nu_{k},\phi_{k})=\sum_{m=1}^{W_{j}}\sum_{n=1}^{W_{j}}a_{m,n}\nu_{k}^{m}\phi_{k}^{n}
\end{equation*}
for $1<j\leq M,1\leq i\leq I_{k,j},1\leq k\leq p.$ Here $1\leq
W_{j}\leq W$. Let
\begin{equation*}
\Omega^{p+1}=\{ \Omega_{l}^{p+1},1\leq l\leq L\}.
\end{equation*}
We define the analytic map $M_{l}^{p+1}$from the master square
$S=(-1,1)^{2}$ to $\Omega_{l}^{p+1}$ and let
\begin{equation*}
u_{l}^{p+1}(M_{l}^{p+1}(\xi,\eta))=\sum_{m=1}^{W}\sum_{n=1}^{W}a_{m,n}\xi^{m}\eta^{n}.
\end{equation*}
Let
$f_{l}^{p+1}(\xi,\eta)=f(X_{l}^{p+1}(\xi,\eta),Y_{l}^{p+1}(\xi,\eta))$
for $1\leq l\leq L$ and $J_{l}^{p+1}(\xi,\eta)$ denote the
Jacobian of the mapping $M_{l}^{p+1}.$ Define $F_{l}^{p+1}
(\xi,\eta)=f_{l}^{p+1}(\xi,\eta)\,\sqrt{J_{l}^{p+1}(\xi,\eta)}$
and let $\hat{F}_{l}^{p+1}(\xi,\eta)$ denote the unique polynomial
which is the orthogonal projection of $F_{l}^{p+1}(\xi,\eta)$ into
the space of polynomials of degree $2W$ in $\xi$ and $\eta$ with
respect to the usual inner product in $H^{2}(S)$.

Next, let the vertex $A_{k}=(x_{k},y_{k}).$ As defined in \S 2 we
have the following relationship between $(\tau_{k},\theta_{k})$:
and $(\nu_{k},\phi_{k})$ coordinates:
\begin{align*}
 \nu_{k} &=\tau_{k},\\[.4pc]
\theta_{k} &
=\frac{1}{(\psi_{u}^{k}-\psi_{l}^{k})}[(\phi_{k}-\psi_{l}^{k})f_{1}^{k}(\hbox{e}^{\nu_{k}})-(\phi_{k}-\psi_{u}^{k})f_{0}^{k}(\hbox{e}^{\nu_{k}})].
\end{align*}
Define $f^{k} (\tau_{k},\theta_{k})=
\hbox{e}^{2\tau_{k}}f(x_{k}+\hbox{e}^{\tau_{k}}\cos\theta_{k},y_{k}+\hbox{e}^{\tau_{k}}\sin\theta_{k})$
for $1\leq k\leq p,$ and
$F_{i,j}^{k}(\nu_{k},\phi_{k})=f^{k}(\tau_{k},\theta_{k})$ for
$(\nu_{k},\phi_{k})\in\hat{\Omega}_{i,j}^{k}.$ Let
$\hat{F}_{i,j}^{k}(\nu_{k},\phi_{k})$ denote the polynomial of
degree $2W_{j}$ in $\nu_{k}$ and $\phi_{k}$ variables which is the
orthogonal projection of $F_{i,j}^{k}(\nu_{k},\phi_{k})$ into the
space of polynomials of degree $2W_{j}$ in $\nu_{k}$ and
$\phi_{k}$ variables with respect to the usual inner product in
$H^{2}(\hat{\Omega}_{i,j}^{k}).$ Here $2\leq j\leq M$.

We now consider the boundary condition $u=g_{k}$ on $\Gamma_{k}$
for $k\in\mathcal{D}$ and let $\big(\frac{\partial u}{\partial
N}\big)_{A}=g_{k}$ on $\Gamma_{k}$ for $k\in\mathcal{N}$. Define
\begin{align*}
\hskip -4pc l_{1}^{k}(\nu_{k})=\begin{cases}
u=g_{k}(x_{k}+\hbox{e}^{\nu_{k}}\cos(f_{_{0}}^{k}(\hbox{e}^{\nu_{k}}))\,,\,
y_{k}+\hbox{e}^{\nu_{k}}\sin(f_{_{0}}^{k}(\hbox{e}^{\nu_{k}}))), \
\  \textrm{for}\,\,
k\in\mathcal{D},\\[.4pc]
\left(\dfrac{\partial u}{\partial
n}\right)_{\tilde{A}^{k}}=\hbox{e}^{\nu_{k}}g_{k}(x_{k}+\hbox{e}^{\nu_{k}}\cos(f_{_{0}}^{k}(\hbox{e}^{\nu_{k}}))\,,\,
y_{k}+e^{\nu_{k}}\sin(f_{_{0}}^{k}(\hbox{e}^{\nu_{k}}))),\\[.7pc]
\quad\, \textrm{for}\,\, k\in\mathcal{N}.
\end{cases}
\end{align*}
Let $\hat{l}$$_{1,j}^{k}(\nu_{k})$ be the orthogonal projection of
$l_{1}^{k}(\nu_{k})$ into the space of polynomials of degree
$2W_{j}$ with respect to the usual inner product on
$H^{2}(\nu_{j}^{k},\nu_{j+1}^{k})$ for $2\leq j\leq M.$

Consider the boundary condition $u=g_{k}$ on
$\Gamma_{k}\cap\partial\Omega^{k-1}$. Define
\begin{equation*}
\hskip -4pc l_{2}^{k}(\nu_{k-1})=\begin{cases}
 u=g_{k}(x_{_{k-1}}+\hbox{e}^{\nu_{k-1}}\cos(f_{_{1}}^{k-1}(\hbox{e}^{\nu_{k-1}}))\,,\, y_{_{k-1}}+\hbox{e}^{\nu_{k-1}}\sin(f_{_{1}}^{k-1}(\hbox{e}^{\nu_{k-1}}))),\\[.4pc]
\quad\, \textrm{for}\,\,  k\in\mathcal{D},\\[.4pc]
\left(\dfrac{\partial u}{\partial
n}\right)_{\tilde{A}^{k}}=\hbox{e}^{\nu_{k-1}}g_{k}(x_{_{k-1}}+\hbox{e}^{\nu_{k-1}}\cos(f_{_{1}}^{k-1}(\hbox{e}^{\nu_{k-1}}))\,,\,
y_{_{k-1}}\\[.8pc]
\quad\,
+\hbox{e}^{\nu_{k-1}}\sin(f_{_{1}}^{k-1}(\hbox{e}^{\nu_{k-1}})))\quad\textrm{for}\,\,
k\in\mathcal{N}.
\end{cases}
\end{equation*}
Let $a_{k}=u(A_{k})$ if $\gamma_{k}$ or
$\gamma_{k+1}\in\mathcal{D}$. We define
$\hat{l}_{2,j}^{k}(\nu_{k-1})$ to be the orthogonal projection of
$l_{2}^{k}(\nu_{k-1})$ into the space of polynomials of degree
$2W_{j}$ with respect to the usual inner product on
$H^{2}(\nu_{j}^{k-1},\nu_{j+1}^{k-1})$ for $2\leq j\leq M.$

Finally, let $\Gamma_{k}\bigcap\partial\Omega_{t}^{p+1}=C_{t}^{k}$
be the image of the mapping $M_{t}^{p+1}$ of $\bar{S}$ onto
$\overline{\Omega}_{t}^{p+1}$ corresponding to the side $\xi=-1.$
Let $o_{t}^{k}(\eta)
=g_{k}(X_{t}^{p+1}(-1,\eta),Y_{t}^{p+1}(-1,\eta)),$ where
$-1\leq\eta\leq1.$ Define $\hat{o}_{t}^{k}(\eta)$ to be the
polynomial of degree $2W$ which is the orthogonal projection of
$o_{t}^{k}(\eta)$ with respect to the usual inner product in
$H^{2}(-1,1).$

Now we formulate the numerical scheme for problems with mixed
boundary conditions.

Let $\{\{ v_{i,j}^{k}(\nu_{k},\phi_{k})\} _{i,j,k},\{
v_{l}^{_{^{p+1}}}(\xi,\eta)\}_{l}\} \in\Pi^{M,W},$ the space of
spectral element functions. Define the functional
\begin{align}
& \mathfrak{\mathcal{\mathfrak{r}}}_{_{\rm vertices}}^{^{M,W}}(\{ v_{i,j}^{k}(\nu_{k},\phi_{k})\}_{i,j,k},\{
v_{l}^{_{^{p+1}}}(\xi,\eta)\}_{l})\nonumber\\[.6pc]
 &\quad\, = \sum_{k=1}^{p}\sum_{j=2}^{M}\sum_{i=1}^{I_{k}}(\rho\mu_{k}^{M+1-j})^{-2\lambda_{k}}\Vert (\mathcal{L}_{i,j}^{k})^{^{a}}v_{i,j}^{k}(\nu_{k},\phi_{k})-\hat{F}_{i,j}^{k}(\nu_{k},\phi_{k})\Vert_{_{0,\hat{\Omega}_{i,j}^{k}}}^{2}\nonumber\\[.6pc]
 &\qquad\, + \sum_{k=1}^{p}\sum_{\gamma_{s}\subseteq\Omega^{k}\cup B_{\rho}^{k},\mu(\hat{\gamma}_{s})<\infty} d(A_{k},\gamma_{s})^{-2\lambda_{k}}\nonumber\\[.6pc]
 &\qquad\, \times (\Vert [v^{k}]\Vert_{_{0,\hat{\gamma}_{s}}}^{2}+\Vert [(v_{\nu_{k}}^{k})^{a}]\Vert_{_{1/2,\hat{\gamma}_{s}}}^{2}+\Vert [(v_{\phi_{k}}^{k})^{a}]\Vert_{_{1/2,\hat{\gamma}_{s}}}^{2})\nonumber\\[.6pc]
 &\qquad\, + \sum_{m\in\mathcal{D}}\sum_{k=m-1}^{m}\sum_{\gamma_{s}\subseteq\partial\Omega^{k}\cap\Gamma_{m},\mu(\hat{\gamma}_{s})<\infty} d(A_{k},\gamma_{s})^{-2\lambda_{k}}(\Vert (v^{k}-h_{k})\nonumber\\[.6pc]
 &\qquad\, -(\hat{l}_{m-k+1}^{m}-a_{k})\Vert_{_{0,\hat{\gamma}_{s}}}^{2} +\Vert v_{\nu_{k}}^{k}-(\hat{l}_{m-k+1}^{m})_{\nu_{k}}\Vert_{_{1/2,\hat{\gamma}_{s}}}^{2})\nonumber \\[.6pc]
 &\qquad +  \sum_{m\in\mathcal{D}}\sum_{k=m-1}^{m}(h_{k}-a_{k})^{2}+\sum_{m\in\mathcal{N}}\sum_{k=m-1}^{m}\nonumber\\[.4pc]
 &\qquad\, \times \sum_{\gamma_{s}\subseteq\partial\Omega^{k}\cap\Gamma_{m},\mu(\hat{\gamma}_{s})<\infty} d(A_{k},\gamma_{s})^{-2\lambda_{k}}\Vert \left(\frac{\partial v^{k}}{\partial n}\right)_{\tilde{A}^{k}}^{a}-\hat{l}_{m-k+1}^{m}\Vert_{_{1/2,\hat{\gamma}_{s}}}^{2}.\label{4.1}
 \end{align}

In the above $\mu(\hat{\gamma}_{s})$ denotes the measure of
$\hat{\gamma}_{s}.$

Next, define
\begin{align}
 & \mathfrak{\mathcal{\mathfrak{r}}}_{_{\rm interior}}^{^{M,W}}(\{ v_{i,j}^{k}(\nu_{k},\phi_{k})\}_{i,j,k},\{
 v_{l}^{_{^{p+1}}}(\xi,\eta)\}_{l})\nonumber\\[.6pc]
 &\quad\, = \sum_{l=1}^{L}\Vert (\mathcal{L}_{l}^{p+1})^{^{a}}v_{l}^{p+1}(\xi,\eta)-\hat{F}_{l}^{p+1}(\xi,\eta)\Vert_{_{0,S}}^{2}\nonumber\\[.6pc]
 & \qquad\,+ \sum_{\gamma_{s}\subseteq\Omega^{^{p+1}}}(\Vert [v^{p+1}]\Vert_{_{0,\gamma_{s}}}^{2}+\Vert [(v_{x_{1}}^{p+1})^{a}]\Vert_{_{1/2,\gamma_{s}}}^{2}+\Vert [(v_{x_{2}}^{p+1})^{a}]\Vert_{_{1/2,\gamma_{s}}}^{2})\nonumber\\[.6pc]
 &\qquad\, + \sum_{l\in\mathcal{D}}\sum_{\gamma_{s}\subseteq\partial\Omega^{^{p+1}}\cap\Gamma_{l}}\left(\!\Vert v^{p+1}-\hat{o}^{l}\Vert_{_{0,\gamma_{s}}}^{2}+\!\left\Vert \left(\frac{\partial v^{p+1}}{\partial T}\right)^{a}-\left(\frac{\partial\hat{o}^{l}}{\partial T}\right)^{a}\right\Vert_{_{1/2,\gamma_{s}}}^{2}\right)\nonumber\\[.6pc]
 &\qquad + \sum_{l\in\mathcal{N}}\sum_{\gamma_{s}\subseteq\partial\Omega^{^{p+1}}\cap\Gamma_{l}}\left\Vert \left(\frac{\partial v^{p+1}}{\partial N}\right)_{A}^{a}-\hat{o}^{l}\right\Vert_{_{1/2,\gamma_{s}}}^{2}.\label{4.1}
\end{align}
Let
\begin{align}
&\mathfrak{\mathfrak{\mathcal{\mathfrak{r}}}}^{^{M,W}}(\{ v_{i,j}^{k}(\nu_{k},\phi_{k})\}_{i,j,k},\{
v_{l}^{_{^{p+1}}}(\xi,\eta)\}_{l})\nonumber\\[.3pc]
&\quad\, = \mathfrak{\mathcal{\mathfrak{r}}}_{_{\rm
vertices}}^{^{M,W}}(\{ v_{i,j}^{k}(\nu_{k},\phi_{k})\}_{i,j,k},\{
v_{l}^{_{^{p+1}}}(\xi,\eta)\}_{l})\nonumber\\[.3pc]
&\qquad\, +\mathcal{\mathfrak{r}}_{_{\rm interior}}^{^{M,W}}(\{
v_{i,j}^{k}(\nu_{k},\phi_{k})\}_{i,j,k},\{
v_{l}^{_{^{p+1}}}(\xi,\eta)\}_{l}).\label{4.3}
\end{align}
We choose as our approximate solution the unique $\{\{
z_{i,j}^{k}(\nu_{k},\phi_{k})\}_{i,j,k},\{
z_{l}^{_{^{p+1}}}(\xi,\eta)\}_{l}\} \in\Pi^{M,W},$ the space of
spectral element functions, which minimizes the functional
$\mathfrak{\mathcal{\mathcal{\mathfrak{r}}}}^{^{M,W}}(\{
v_{i,j}^{k}(\nu_{k},\phi_{k})\}_{i,j,k},\{
v_{l}^{_{^{p+1}}}(\xi,\eta)\}_{l})$ over all $\{ \{
v_{i,j}^{k}(\nu_{k},\phi_{k})\}_{i,j,k},\{
v_{l}^{_{^{p+1}}}(\xi,\eta)\}_{l}\} $.

A brief description of the solution procedure is now given; a more
detailed examination is provided in \S6. The above method is
essentially a least-squares method and the solution can be
obtained by using preconditioned conjugate gradient techniques
(PCGM) to solve the normal equations. To be able to do so we must
be able to compute the residuals in the normal equations
inexpensively. In \cite{prsb,tomarth} it has been shown how to
compute these efficiently on a distributed memory parallel
computer, without having to filter the coefficients of the
differential operator and the data. The evaluation of the
residuals on each element requires the interchange of boundary
values between neighbouring elements.

The values of the spectral element functions at the vertices of
the polygonal domain constitute the set of common boundary values
$U_{B}.$ Since the dimension of the set of common boundary values
is so small a nearly exact approximation to the Schur Complement
matrix can be computed. Now on the subspace of spectral element
functions which vanish at the set of common boundary values it is
possible to define a preconditioner for the matrix in the normal
equations such that the condition number of the preconditioned
system is $O((\ln W)^{2}).$ Moreover, the preconditioner is a
block diagonal matrix such that each diagonal block corresponds to
a different element, and so can be easily inverted.

Hence an exponentially accurate approximation $\mathbb{S}^{a}$ to
the Schur Complement matrix $\mathbb{S}$ can be computed using
$O(W\ln W)$ iterations of the PCGM. To solve the normal equations
the residual in the equations for the Schur Complement
$\mathbb{S}U_{B}=h_{B}$ must be computed to exponential accuracy
and this can be done using $O(W\ln W)$ iterations of the PCGM. The
common boundary values $U_{B}$ are then given by
$U_{B}=(\mathbb{S}^{a})^{-1}h_{B}.$ The remaining values can then
be obtained using $O(W\ln W)$ iterations of the PCGM.

\section{Error estimates}

\setcounter{equation}{0}

\setcounter{defin}{0}

Let $\{\{z_{i,j}^{k}(\nu_{k},\phi_{k})\}_{i,j,k},
\{z_{l}^{_{^{p+1}}}(\xi,\eta)\}_{l}\}$ minimize
$\mathfrak{\mathcal{\mathcal{\mathfrak{r}}}}^{^{M,W}}
(\{v_{i,j}^{k}(\nu_{k},\phi_{k})\}_{i,j,k},\{v_{l}^{_{^{p+1}}}(\xi,\eta)\}_{l})$
over all $\{\{v_{i,j}^{k}(\nu_{k},\phi_{k})\}_{i,j,k},
\{v_{l}^{_{^{p+1}}}(\xi,\eta)\}_{l}\} \in\Pi^{M,W}$, the space of
spectral element functions. Here $z_{i,1}^{k}=b_{k}$ for all $i,$
$z_{i,j}^{k}(\nu_{k},\phi_{k})$ is a polynomial in $\nu_{k}$ and
$\phi_{k}$ of degree $W_{j}$, $W_{j}\leq W$ and
$z_{l}^{_{^{p+1}}}(\xi,\eta)$ is a polynomial in $\xi$ and $\eta$
of degree $W$ as defined in \S3. We choose $W$ proportional to
$M.$ Then we have the following error estimate.

\begin{theo}[\!]
Let $a_{k}=u(A_{k})$. Let $U_{i,j}^{k}(\nu_{k},\phi_{k}) =
u(\nu_{k},\phi_{k})$ for $(\nu_{k},\phi_{k})\in
\hat{\Omega}_{i,j}^{k}$ and $U_{l}^{p+1}(\xi,\eta)=u(\xi,\eta)$
for $(\xi,\eta)\in S$. Let $\alpha j\leq W_{j}\leq W$ for some
positive $\alpha$ for $j>2$. Then there exists positive constants
$C$ and $b$ such that for $W$ large enough the estimate
\begin{align*}
&\sum_{k=1}^{p}|b_{k}-a_{k}|^{2}+\sum_{k=1}^{p}\sum_{j=2}^{M}
\sum_{i=1}^{I_{k}}(\rho\mu_{k}^{M+1-j})^{-2\lambda_{k}}
\end{align*}
\begin{align}
&\quad\, \times \Vert
(z_{i,j}^{k}-U_{i,j}^{k})(\nu_{k},\phi_{k})-(b_{k}-a_{k})
\Vert_{2,\hat{\Omega}_{i,j}^{k}}^{2}\nonumber\\[.5pc]
&\quad\, + \sum_{l=1}^{L}\Vert (z_{l}^{p+1}-U_{l}^{p+1})
(\xi,\eta)\Vert_{2,S}^{2}\leq C\, {\rm e}^{-bW}\label{5.1}
\end{align}
holds.
\end{theo}
We use the differentiability estimates stated in Proposition~2.1
to prove the result. The proof of the above Theorem is very
similar to the proof of Theorem~3.1 in \cite{tomarth} and hence is
omitted.\hfill $\blacksquare$

\begin{remarr}
We can construct a set of corrections
$\{\{c_{i,j}^{k}(\nu_{k},\phi_{k})\}_{i,j,k},
\{c_{l}^{_{^{p+1}}}(\xi,\eta)\}_{l}\} \in\Pi^{M,W},$ the set of
spectral element functions, so that corrected solution
$\{\{\hat{z}_{i,j}^{k}(\nu_{k},\phi_{k})\}_{i,j,k}$,
$\{\hat{z}_{l}^{_{^{p+1}}}(\xi,\eta)\}_{l}\}$ defined by
\begin{align*}
&\{\{\hat{z}_{i,j}^{k}(\nu_{k},\phi_{k})\}_{i,j,k},
\{\hat{z}_{l}^{_{^{p+1}}}(\xi,\eta)\}_{l}\}\\[.3pc]
&\quad\, = \{\{z_{i,j}^{k} (\nu_{k},\phi_{k})\}_{i,j,k},
\{z_{l}^{_{^{p+1}}} (\xi,\eta)\}_{l}\} +\{\{c_{i,j}^{k}(\nu_{k},
\phi_{k})\}_{i,j,k},\{c_{l}^{_{^{p+1}}}(\xi,\eta)\}_{l}\}
\end{align*}
is conforming and belongs to $H^{1}(\Omega).$ These corrections
are defined in \S3.5 of \cite{tomarth}. Then the error estimate
\begin{equation*}
\Vert (u-\hat{z})(x,y)\Vert_{1,\Omega}\leq C\, {\rm e}^{-bW}
\end{equation*}
holds for $W$ large enough. Here $C$ and $b$ denote constants.
These constructions are similar to Lemma~4.57 in \cite{schwab}.
\end{remarr}

\section{Parallelization and preconditioning}

\setcounter{equation}{0}

\setcounter{defin}{0}

Let $U$ be a vector assembled from $\{g_{k}\}_{k=1}^{p},$ where
$u_{i,1}^{k} =g_{k}$ for all $i,$ and the values of
$\{\{u_{i,j}^{k}(\nu_{k},\phi_{k})\}_{i,j,k},
\{u_{l}^{p+1}(\xi,\eta)\}_{_{l}}\} $ at the
Gauss--Lobatto--Legendre points are arranged in lexicographic
order for $1\leq k\leq p,$ $2\leq j\leq J_{k}$, $1\leq i\leq
I_{k,j}.$ Let $\{\{z_{i,j}^{k}(\nu_{k},\phi_{k})\}_{i,j,k},
\{z_{l}^{_{^{p+1}}}(\xi,\eta)\}_{l}\}$ minimize
$\mathfrak{r}^{^{M,W}}(\{v_{i,j}^{k}(\nu_{k},\phi_{k})\}_{i,j,k} ,
\{v_{l}^{_{^{p+1}}}(\xi,\eta)\}_{l})$ over all
$\{\{v_{i,j}^{k}(\nu_{k},\phi_{k})\}_{i,j,k},
\{v_{l}^{_{^{p+1}}}(\xi,\eta)\}_{l}\} \in\Pi^{M,W}$, the space of
spectral element functions.

Let $U_{B}$ denote the values $\{g_{k}\}_{k=1}^{p}$ and $U_{I}$
the remaining values of $U$. We now define a quadratic form
\begin{align}
&\mathcal{Z}^{^{M,W}}(\{u_{i,j}^{k}(\nu_{k},\phi_{k})\}_{i,j,k},
\{u_{l}^{_{^{p+1}}}(\xi,\eta)\}_{l})\nonumber\\[.3pc]
&\quad\, = \sum_{k=1}^{p}|g_{k}|^{2}+\sum_{k=1}^{p}\sum_{j=2}^{M}
\sum_{i=1}^{I_{k}}(\rho\mu_{k}^{M+1-j})^{-2\lambda_{k}}\Vert
u_{i,j}^{k}(\xi,\eta)-g_{k}\Vert_{2,S}^{2}\nonumber\\[.3pc]
&\qquad\, +\sum_{l=1}^{L}\Vert
u_{l}^{p+1}(\xi,\eta)\Vert_{2,S}^{2}.\label{5.1}
\end{align}
It should be noted that $u_{i,1}^{k}(\nu_{k},\phi_{k})=g_{k}$ for
$1\leq i\leq I_{k}.$ Moreover for $j\leq M$, $\xi$ is a linear
function of $\nu_{k}$ and $\eta$ is a linear function of
$\phi_{k}$ such that the linear mapping $M_{i,j}^{k}(\xi,\eta)$
maps the master square $S$ onto
$\hat{\Omega}_{i,j}^{k}.$\pagebreak

To solve the minimization problem we have to solve a system of
equations of the form
\begin{align}
A Z=h.\label{6.2}
\end{align}

Here $A$ is a symmetric positive definite matrix and
\begin{align}\label{5.3}
\mathcal{V}^{^{M,W}}(\{u_{i,j}^{k}(\nu_{k},\phi_{k})\}_{i,j,k},
\{u_{l}^{_{^{p+1}}}(\xi,\eta)\}_{l})=U^{T}A\, U,\hskip
-1pc\phantom{0}
\end{align}
where $\mathcal{V}^{^{M,W}}(\{u_{i,j}^{k}(\nu_{k},
\phi_{k})\}_{i,j,k}, \{u_{l}^{_{^{p+1}}}(\xi,\eta)\}_{l})$ is as
defined in (3.15) in \S3.

Now $A$ has the form
\begin{align}\label{5.4}
A=\left[\begin{array}{cc}
A_{II} & A_{IB}\\[.2pc]
A_{BI} & A_{BB}\end{array}\right]\hskip -1pc\phantom{0}
\end{align}
corresponding to the decomposition of $U$ as
\begin{equation*}
U=\left[\begin{array}{c}
U_{I}\\[.2pc]
U_{B}\end{array}\right],
\end{equation*}
and $h$ has the form
\begin{align*}
h=\left[\begin{array}{c}
h_{I}\\[.2pc]
h_{B}\end{array}\right].
\end{align*}
To solve the matrix equation (6.2) we use the block L-U
factorization of $A$, viz.
\begin{align}\label{5.5}
A=\left[\begin{array}{cc}
I &0\\[.2pc]
A_{IB}^{T}A_{II}^{-1}
&I\end{array}\right]\,\left[\begin{array}{cc}
A_{II} &0\\[.2pc]
0 &\mathbb{S}\end{array}\right]\,\left[\begin{array}{cc}
I &A_{II}^{-1}A_{IB}\\[.2pc]
0 &I\end{array}\right],\hskip -1pc\phantom{0}
\end{align}
where the Schur Complement matrix $\mathbb{S}$ is defined as
\begin{align}
\mathbb{S}=A_{BB}-A_{IB}^{T}A_{II}^{-1}A_{IB}.\label{5.6}
\end{align}

To solve the matrix equation (6.2) based on the L-U factorization
of $A$ given in (6.5) reduces to solving the system of equations
\begin{align}\label{5.7}
\mathbb{S} Z_{B}=\tilde{h}_{B},\hskip -1pc\phantom{0}
\end{align}
where
\begin{align}
\tilde{h}_{B}=h_{B}-A_{IB}^{T}A_{II}^{-1}h_{I}.\label{6.8}
\end{align}
The feasibility of such a process depends on our being able to
compute $A_{IB}V_{B}$, $A_{II}V_{I}$ and $A_{BB}V_{B}$ for any
$V_{I},V_{B}$ efficiently and this can always be done since $AV$
can be computed inexpensively as explained in ch.~3 of
\cite{tomarth}.

However in addition to this it is imperative that we should be
able to construct effective preconditioners for the matrix
$A_{II}$ so that the condition number of the preconditioned system
is as small as possible. If this can be done then it will be
possible to compute $A_{II}^{-1}V_{I}$ efficiently using the
preconditioned conjugate gradient method (PCGM) for any vector
$V_{I}$.

Consider the space of spectral element functions $\Pi_{0}^{M,W},$
such that for $\{\{u_{i,j}^{k}(\nu_{k},\phi_{k})\}_{i,j,k},$
$\{u_{l}^{_{^{p+1}}}(\xi,\eta)\}_{l}\} \in\Pi_{0}^{M,W}$ we have
$u_{i,1}^{k}=0$ for all $i$ and $k.$ Let $U$ be the vector
corresponding to the spectral element function
$\{\{u_{i,j}^{k}(\nu_{k},\phi_{k})\}_{^{i,j,k}},\{u_{l}^{_{^{p+1}}}(\xi,\eta)\}_{l}\}$.
Then $U_{B}=0$ and $U=\Big[\begin{smallmatrix}
U_{I}\\[.1pc]
0\end{smallmatrix}\Big]$ and so
\begin{align}
\mathcal{V}^{^{M,W}}(\{u_{i,j}^{k}(\nu_{k},\phi_{k})\}_{i,j,k},
\{u_{l}^{_{^{p+1}}}(\xi,\eta)\}_{l})=U_{I}^{T}A_{II}\,
U_{I}.\label{5.9}
\end{align}
Now using Theorem~3.1 we have the following result.

Let $\{\{u_{i,j}^{k}(\nu_{k},\phi_{k})\}_{i,j,k},
\{u_{l}^{_{^{p+1}}}(\xi,\eta)\}_{l}\} \in\Pi_{0}^{M,W}$. Then the
estimate
\begin{align}
&\sum_{k=1}^{p}\sum_{j=2}^{M}\sum_{i=1}^{I_{k}}
(\rho\mu_{k}^{M+1-j})^{-2\lambda_{k}} \Vert
u_{i,j}^{k}(\xi,\eta)\Vert_{2,S}^{2}+\sum_{l=1}^{L}\Vert
u_{l}^{p+1}(\xi,\eta)\Vert_{2,S}^{2}\nonumber\\[.5pc]
&\quad\, \leq C(\ln
W)^{2}\,\mathcal{V}^{^{M,W}}(\{u_{i,j}^{k}(\nu_{k},
\phi_{k})\}_{i,j,k},\{u_{l}^{_{^{p+1}}}(\xi,\eta)\}_{l})\label{6.10}
\end{align}
holds for $W$ large enough. In the above, $u_{i,1}^{k}=0$ for
$1\leq k\leq p$ and $1\leq i\leq I_{k}$.

Let us define the quadratic form
\begin{align}
&\mathcal{U}^{M,W}(\{u_{i,j}^{k}(\nu_{k},\phi_{k})\}_{i,j,k},
\{u_{l}^{_{^{p+1}}}(\xi,\eta)\}_{l})\nonumber\\[.5pc]
&\quad\, =\sum_{k=1}^{p}\sum_{j=2}^{M}\sum_{i=1}^{I_{k}}
(\rho\mu_{k}^{M+1-j})^{-2\lambda_{k}}\Vert
u_{i,j}^{k}(\xi,\eta)\Vert_{2,S}^{2}+\sum_{l=1}^{L}\Vert
u_{l}^{p+1}(\xi,\eta)\Vert_{2,S}^{2}\label{6.11}
\end{align}
for all $\{\{u_{i,j}^{k}(\nu_{k},\phi_{k})\}_{i,j,k},
\{u_{l}^{_{^{p+1}}}(\xi,\eta)\}_{l}\} \in\Pi_{0}^{M,W}$.

Now using the trace theorems for Sobolev spaces it can be
concluded that there exists a constant $K$ such that
\begin{align}
&\mathcal{V}^{^{M,W}}(\{u_{i,j}^{k}(\nu_{k},\phi_{k})\}_{i,j,k},
\{u_{l}^{_{^{p+1}}}(\xi,\eta)\}_{l})\nonumber\\[.5pc]
&\quad\, \leq
K\,\mathcal{U}^{M,W}(\{u_{i,j}^{k}(\nu_{k},\phi_{k})\}_{i,j,k},
\{u_{l}^{_{^{p+1}}}(\xi,\eta)\}_{l})\label{6.12}
\end{align}
for $\{\{u_{i,j}^{k}(\nu_{k},\phi_{k})\}_{i,j,k},
\{u_{l}^{_{^{p+1}}}(\xi,\eta)\}_{l}\} \in\Pi_{0}^{M,W}$.

Hence using (6.10) and (6.12) it follows that there exists a
constant $C$ such that
\begin{align}
&\frac{1}{C}
\mathcal{V}^{^{M,W}}(\{u_{i,j}^{k}(\nu_{k},\phi_{k})\}_{i,j,k},
\{u_{l}^{_{^{p+1}}}(\xi,\eta)\}_{l})\nonumber\\[.5pc]
&\quad\, \leq
\mathcal{U}^{M,W}(\{u_{i,j}^{k}(\nu_{k},\phi_{k})\}_{i,j,k},
\{u_{l}^{_{^{p+1}}}(\xi,\eta)\}_{l})\nonumber\\[.5pc]
&\quad\, \leq C\,(\ln
W)^{2}\,\mathcal{V}^{^{M,W}}(\{u_{i,j}^{k}(\nu_{k},
\phi_{k})\}_{i,j,k},\{u_{l}^{_{^{p+1}}}(\xi,\eta)\}_{l})\label{5.13}
\end{align}
for all $\{\{u_{i,j}^{k}(\nu_{k},\phi_{k})\}_{i,j,k},
\{u_{l}^{_{^{p+1}}}(\xi,\eta)\}_{l}\} \in\Pi_{0}^{M,W}$.

Thus the two forms
$\mathcal{V}^{^{M,W}}(\{u_{i,j}^{k}(\nu_{k},\phi_{k})\}_{i,j,k},
\{u_{l}^{_{^{p+1}}}(\xi,\eta)\}_{l})$ and
$\mathcal{U}^{^{M,W}}(\{u_{i,j}^{k}(\nu_{k},\phi_{k})\}_{i,j,k},$
$\{u_{l}^{_{^{p+1}}}(\xi,\eta)\}_{l})$ are spectrally equivalent.

We can now use the quadratic form
$\mathcal{U}^{^{M,W}}(\{u_{i,j}^{k}(\nu_{k},\phi_{k})\}_{i,j,k},
\{u_{l}^{_{^{p+1}}}(\xi,\eta)\}_{l})$ which consists of a
decoupled set of quadratic forms on each element as a
preconditioner for $A_{II}$. This can be done by inverting the
block diagonal matrix representation for
$\mathcal{U}^{^{M,W}}(\{u_{i,j}^{k}(\nu_{k},
\phi_{k})\}_{i,j,k},\{u_{l}^{_{^{p+1}}}(\xi,\eta)\}_{l})$.

Now from (6.13) we can conclude that if we were to compute
$(A_{II})^{-1}U_{I}$ using the PCGM then the condition number of
the preconditioned matrix would be $O((\ln W)^{2})$. Hence, to
compute $(A_{II})^{-1}U_{I}$ to an accuracy of $O({\rm e}^{-bW})$
would require $O(W\, \ln W)$ iterations of the PCGM.

We now return to the steps involved in solving the system of
equations~(6.2). As a first step it would be necessary to solve
the much smaller system of equations~(6.7). Here the dimension of
the vector $Z_{B}$ is $p$, the number of vertices of the domain
$\Omega$. Now to be able to solve (6.7) to an accuracy of $O({\rm
e}^{-bW})$ using PCGM the residual
\begin{equation*}
R_{B}=\mathbb{S}U_{B}-\tilde{h}_{B}
\end{equation*}
needs to be computed with the same accuracy and in an efficient
manner. The bottleneck in computing $R_{B}$ consists in computing
$(A_{II})^{-1}A_{IB}U_{B}$ to an accuracy of $O({\rm e}^{-bW})$
and it has already been seen that this can be done using $O(W\,
\ln W)$ iterations of the PCGM for computing
$(A_{II})^{-1}A_{IB}U_{B}$ for a given vector $U_{B}.$

We now show that it is possible to explicitly construct the Schur
Complement matrix $\mathbb{S}$ in $O(W\,\ln W)$ iterations of the
PCGM. $\mathbb{S}$ is a $p\times p$ matrix. Let $e_{k}$ be a
column vector of dimension $p$ with 1 in the $k$th place and 0
elsewhere. Let $\mathbb{S}_{k}=\mathbb{S}e_{k}.$

Then the Schur Complement matrix $\mathbb{S}$ can be written as
\begin{align*}
\mathbb{S}=[\mathbb{S}_{1},\mathbb{S}_{2},\dots,\mathbb{S}_{p}].
\end{align*}
Now by a well known result on the Schur Complement we have
\begin{align*}
U_{B}^{T}\mathbb{S}U_{B} &= \min_{_{V\hbox{\rm :}\ V_{B}=U_{B}}}V^{T}AV\\[.3pc]
&= \min_{v_{i,j}^{k}\hbox{\rm :}\ v_{i,1}^{k} =
g_{k}}\mathcal{V}^{^{M,W}}(\{v_{i,j}^{k}(\nu_{k},\phi_{k})\}_{i,j,k},
\{v_{l}^{_{^{p+1}}}(\xi,\eta)\}_{l}).
\end{align*}
Here $U_{B}=[g_{1},g_{2},\dots,g_{p}]^{T}.$ Hence using
Theorem~3.1 we conclude that
\begin{equation*}
U_{B}^{T}\mathbb{S}U_{B}\geq\frac{C}{(\ln W)^{2}}\Vert
U_{B}\Vert^{2}.
\end{equation*}
And so we obtain
\begin{align}
\Vert \mathbb{S}^{-1}\Vert \leq C(\ln W)^{2}.\label{5.14}
\end{align}
Here the norm denoted is the matrix norm induced by the Euclidean
norm. Now
\begin{align*}
\mathbb{S}_{k}=\mathbb{S}e_{k}=(A_{BB}-A_{IB}^{T}A_{II}^{-1}A_{IB})\,\,
{\rm e}_{k}.
\end{align*}
Let $(\mathbb{S}_{k})^{a}$ be the approximation to
$\mathbb{S}_{k}$ computed using $O(W\ln W)$ iterations of the PCGM
to compute $A_{II}^{-1}A_{IB}e_{k}.$ Then
\begin{equation*}
\Vert \mathbb{S}_{k}-\mathbb{S}_{k}^{a}\Vert =O({\rm e}^{-bW}).
\end{equation*}
Let $\mathbb{S}^{a}$ denote the matrix
\begin{equation*}
\mathbb{S}^{a}=[\mathbb{S}_{1}^{a},\mathbb{S}_{2}^{a},\dots,
\mathbb{S}_{p}^{a}].
\end{equation*}
Clearly
\begin{equation*}
\Vert \mathbb{S}-\mathbb{S}^{a}\Vert =O({\rm e}^{-bW}).
\end{equation*}
Now to compute $\mathbb{S}^{a}$ requires $O(W\ln W)$ iterations of
the PCGM since $p$ is a fixed constant. Hence we can solve (6.7)
as
\begin{equation*}
\mathbb{S}Z_{B}=\tilde{h}_{B}
\end{equation*}
by replacing $\mathbb{S}$ by the matrix $\mathbb{S}^{a}$. Let
$Z_{B}^{a}$ be the solution of
\begin{equation*}
\mathbb{S}^{a}Z_{B}^{a}=\tilde{h}_{B}.
\end{equation*}
Since
\begin{equation*}
\mathbb{S}^{a}=\mathbb{S}+\delta\mathbb{S},
\end{equation*}
we have
\begin{equation*}
(\mathbb{S}^{a})^{-1}=(I+\mathbb{S}^{-1}\delta\mathbb{S})\mathbb{S}^{-1}.
\end{equation*}
Thus
\begin{equation*}
\Vert \mathbb{S}^{-1}-(\mathbb{S}^{a})^{-1}\Vert \leq2\Vert
\mathbb{S}^{-1}\Vert ^{2}\Vert \delta\mathbb{S}\Vert \leq O((\ln
W)^{4})\Vert \delta\mathbb{S}\Vert
\end{equation*}
for $\Vert \delta\mathbb{S}\Vert$ small enough.

Hence
\begin{equation*}
\Vert \mathbb{S}^{-1}-(\mathbb{S}^{a})^{-1}\Vert =O({\rm
e}^{-bW}).
\end{equation*}
Therefore
\begin{equation*}
\Vert Z_{B}^{a}-Z_{B}\Vert =O({\rm e}^{-bW}).
\end{equation*}
Having solved for $Z_{B}$ we obtain $Z_{I}$ by solving
\begin{equation*}
A_{II}Z_{I}=h_{I}-A_{IB}Z_{B}
\end{equation*}
using $O(W\ln W)$ iterations of the PCGM. Hence the solution $Z$
can be obtained to exponential accuracy using $O(W\ln W)$
iterations of the PCGM.

We shall now briefly examine the complexity of the solution
procedure for the h-p finite element method. Since finite elements
have to be continuous along the sides of the elements, the
cardinality of the set of common boundary value is large in the
h-p finite element method. Let $\mathbb{S}$ denote the Schur
Complement matrix for the h-p finite element method. In
\cite{babcra,guocao} it has been shown that an approximation
$\mathbb{S}^{a}$ to $\mathbb{S}$ can be obtained such that the
condition number $\chi$ of the preconditioned system satisfies
\begin{equation*}
\chi\leq C(1+(\ln W)^{2}),
\end{equation*}
where $C$ denotes a constant. Then to solve $\mathbb{S} U_{B} =
h_{B}$ to an accuracy $O({\rm e}^{-bW})$ will require $O(W\ln W)$
iterations of the PCGM using $\mathbb{S}^{a}$ as a preconditioner.
Now to compute the residual in the Schur Complement system to an
accuracy of $O({\rm e}^{-bW})$ requires $O(W)$ iterations of the
PCGM to compute $A_{II}^{-1}A_{IB}V_{B}.$ Hence we would need to
perform $O(W^{2}\, \ln W)$ iterations of the PCGM for computing
$A_{II}^{-1}V_{I}$, where $V_{I}$ will vary after every sequence
of $O(W\ln W)$ steps. So the h-p finite element method requires
$O(W^{2}\ln W)$ iterations of the PCGM to obtain the solution.

Hence the proposed method is faster than h-p finite element method
by a factor of $O(W)$.

\section{Technical results}

\setcounter{equation}{0}

\setcounter{defin}{0}

\begin{lem}
Let $\{\{u_{i,j}^{k}(\nu_{k},\phi_{k})\}_{i,j,k},
\{u_{l}^{_{^{p+1}}}(\xi,\eta)\}_{l}\} \in\Pi^{M,W}$. Then there
exists $\{\{v_{i,j}^{k}$ $(\nu_{k},\phi_{k})\}_{i,j,k},$
$\{v_{l}^{_{^{p+1}}}(\xi,\eta)\}_{l}\}$ such that $v_{i,1}^{k}=0$
for all $i,k$, $v_{i,j}^{k}\in H^{2}(\hat{\Omega}_{i,j}^{k})$ for
$2\leq j\leq M$ and all $i$ and $k$, $v_{l}^{p+1}\in H^{2}(S)$ for
$l=1,2,\dots,L$ and $w=u+v\in H_{\beta}^{2,2}(\Omega).$ Moreover
the estimate
\begin{align}
&\sum_{k=1}^{p}\sum_{i=1}^{I_{k}}\sum_{j=2}^{M}(\rho\mu_{k}^{M+1-j})^{-2\lambda_{k}}\Vert
v_{i,j}^{k}(\nu_{k},\phi_{k})\Vert_{_{2,\hat{\Omega}_{i,j}^{k}}}^{2}+\sum_{l=1}^{L}\Vert
v_{l}^{p+1}(\xi,\eta)\Vert_{_{2,S}}^{2}\nonumber\\[.3pc]
&\quad\, \leq C(\ln
W)^{2}\Bigg(\sum_{k=1}^{p}\sum_{\gamma_{s}\subseteq\Omega^{k}\cup
B_{\rho}^{k},\mu(\hat{\gamma}_{s})<\infty}
d(A_{k},\gamma_{s})^{-2\lambda_{k}}\nonumber\\[.3pc]
&\qquad\, \times (\Vert
[u^{k}]\Vert_{_{0,\hat{\gamma}_{s}}}^{2}+\Vert
[(u_{\nu_{k}}^{k})^{a}]\Vert_{_{1/2,\hat{\gamma}_{s}}}^{2}+\Vert
[(u_{\phi_{k}}^{k})^{a}]\Vert_{_{1/2,\hat{\gamma}_{s}}}^{2})\nonumber\\[.3pc]
&\qquad\, + \sum_{\gamma_{s}\subseteq\Omega^{^{p+1}}}(\Vert
[u^{p+1}]\Vert_{_{0,\gamma_{s}}}^{2}+\Vert
[(u_{x_{1}}^{p+1})^{a}]\Vert_{_{1/2,\gamma_{s}}}^{2}+\Vert
[(u_{x_{2}}^{p+1})^{a}]\Vert_{_{1/2,\gamma_{s}}}^{2})\Bigg)\nonumber\\[.3pc]
&\qquad\, +
\varepsilon_{_{W}}\,\Bigg(\sum_{k=1}^{p}\sum_{i=1}^{I_{k}}
\sum_{j=2}^{M}(\rho\mu_{k}^{M+1-j})^{-2\lambda_{k}}\Vert
u_{i,j}^{k}(\nu_{k},\phi_{k})-h_{k}\Vert_{_{2,
\hat{\Omega}_{i,j}^{k}}}^{2}\nonumber\\[.3pc]
&\qquad\, +\sum_{l=1}^{L}\Vert
u_{l}^{p+1}(\xi,\eta)\Vert_{_{2,S}}^{2}\Bigg)\label{8.1}
\end{align}
holds. Here $\varepsilon_{_{W}}$ is exponentially small in $W$.
\end{lem}

We first make a correction
$\{\{r_{i,j}^{k}(\nu_{k},\phi_{k})\}_{i,j,k},
\{r_{l}^{_{^{p+1}}}(\xi,\eta)\}_{l}\}$ such that $r_{i,1}^{k}=0$
for all $i$ and $k$ and at the vertices $\hat{Q}_{l}$ for
$l=1,\dots,4$ of $\hat{\Omega}_{i,j}^{k}$,
\begin{subequations}
\begin{align}\label{8.2a}
(u_{i,j}^{k}+r_{i,j}^{k})(\hat{Q}_{l})&=\bar{u}(\hat{Q}_{l}),\nonumber\\[.3pc]
((u_{i,j}^{k})_{\nu_{k}}+(r_{i,j}^{k})_{\nu_{k}})(\hat{Q}_{l})&=\bar{u}_{\nu_{k}}(\hat{Q}_{l}),\nonumber\\[.3pc]
((u_{i,j}^{k})_{\phi_{k}}+(r_{i,j}^{k})_{\phi_{k}})(\hat{Q}_{l})&=\bar{u}_{\phi_{k}}(\hat{Q}_{l}),\hskip
-1pc\phantom{0}
\end{align}
provided $Q_{l}$ is not a vertex of $\Omega_{i,1}^{k}$ for all
$i,k$. If $Q_{l}$ is a vertex of $\Omega_{i,1}^{k}$ choose
$r_{i,2}^{k}$ such that
\begin{align}
(u_{i,2}^{k}+r_{i,2}^{k})(\hat{Q}_{l})&=u_{i,1}^{k}(\hat{Q}_{l}),\nonumber\\[.3pc]
((u_{i,2}^{k})_{\nu_{k}}+(r_{i,2}^{k})_{\nu_{k}})(\hat{Q}_{l})&=(u_{i,1}^{k})_{\nu_{k}}(\hat{Q}_{l}),\nonumber\\[.3pc]
((u_{i,2}^{k})_{\phi_{k}}+(r_{i,2}^{k})_{\phi_{k}})(\hat{Q}_{l})&=(u_{i,1}^{k})_{\phi_{k}}(\hat{Q}_{l}).\label{8.2b}
\end{align}
\end{subequations}
Here $\overline{s}(\hat{Q}_{l})$ denotes the average of the values
of $s$ at $\hat{Q}_{l}$ over all the elements which have
$\hat{Q}_{l}$ as a vertex.

We can find a polynomial $r_{i,j}^{k}(\nu_{k},\phi_{k})$ on
$\hat{\Omega}_{i,j}^{k}$ such that
$r_{i,j}^{k}(\hat{Q}_{l})=a_{l},
(r_{i,j}^{k})_{\nu_{k}}(\hat{Q}_{l})=
b_{l}$,$(r_{i,j}^{k})_{\phi_{k}}(\hat{Q}_{l})=c_{l}$ for
$l=1,\dots,4.$ Here the values $a_{l},b_{l},c_{l}$ are defined by
(7.2). Moreover $r_{i,j}^{k}$ is a polynomial of degree less than
or equal to four and the estimate\pagebreak
\begin{align}\label{8.3}
&\Vert r_{i,j}^{k}(\nu_{k},\phi_{k})
\Vert_{_{2,\hat{\Omega}_{i,j}^{k}}}^{2}\leq
C\,\left(\sum_{l=1}^{4}|a_{l}|^{2}+|b_{l}|^{2}+|c_{l}|^{2}\right)\hskip
-1pc\phantom{0}
\end{align}
holds for $j\geq2$ and all $i$ and $k$. Next consider
$\Omega_{l}^{p+1}\in\Omega^{p+1}.$ Now
\begin{align*}
&(u_{l}^{p+1})_{x_{1}}=(u_{l}^{p+1})_{\xi}\,\xi_{x_{1}}+(u_{l}^{p+1})_{\eta}\,\eta_{x_{1}},\,\,\textrm{and}\\[.5pc]
&(u_{l}^{p+1})_{x_{2}}=(u_{l}^{p+1})_{\xi}\,\xi_{x_{2}}+(u_{l}^{p+1})_{\eta}\,\eta_{x_{2}}.
\end{align*}
Let $\hat{\xi}_{x_{1}}$, $\hat{\xi}_{x_{2}}$, $\hat{\eta}_{x_{1}}$
and $\hat{\eta}_{x_{2}}$ denote the polynomials of degree $W$ in
$\xi$ and $\eta\,$ separately which are the approximations to
$\xi_{x_{1}}$, $\xi_{x_{2}}$, $\eta_{x_{1}}\,$and $\eta_{x_{2}}$
in the space of polynomial of degree $W$ as defined in
Theorem~4.46 of \cite{schwab}.

Let $P_{j}$, $j=1,\dots,4$ denote the vertices of $S.$ Then
$\hat{\xi}_{x_{i}}(P_{j})=\xi_{x_{i}}(P_{j})$ and
$\hat{\eta}_{x_{i}}(P_{j})=\eta_{x_{i}}(P_{j})$ for $i=1,2$ and
$j=1,\dots,4.$ Now
\begin{align*}
&(u_{l}^{p+1})_{x_{1}}^{a}=(u_{l}^{p+1})_{\xi}\hat{\xi}_{x_{1}}+(u_{l}^{p+1})_{\eta}\,\hat{\eta}_{x_{1}},\,\:\textrm{and}\\[.5pc]
&(u_{l}^{p+1})_{x_{2}}^{a}=(u_{l}^{p+1})_{\xi}\,\hat{\xi}_{x_{2}}+(u_{l}^{p+1})_{\eta}\,\hat{\eta}_{x_{2}}.
\end{align*}
Hence
$(u_{l}^{p+1})_{x_{i}}^{a}(P_{j})=(u_{l}^{p+1})_{x_{i}}(P_{j}),$
for $i=1,2$ and $j=1,\dots,4.$ Therefore we can find a polynomial
$r_{l}^{p+1}(\xi,\eta)$ on
$S=(M_{l}^{p+1})^{-1}(\Omega_{l}^{p+1})$ such that for
$j=1,\dots,4$,
\begin{align*}
&(u_{l}^{p+1}+r_{l}^{p+1})(P_{j})=\bar{u}(P_{j}),\\[.5pc]
&((u_{l}^{p+1})_{x_{1}}+(r_{l}^{p+1})_{x_{1}})(P_{j})=\bar{u}_{x_{1}}(P_{j}),\,\,\textrm{and}\\[.5pc]
&((u_{l}^{p+1})_{x_{2}}+(r_{l}^{p+1})_{x_{2}})(P_{j})=\bar{u}_{x_{2}}(P_{j}).
\end{align*}
Now let $q(y)$ be a polynomial of degree $W$ defined on
$I=(-1,1)$. Then by Theorem~4.79 of \cite{schwab}
\begin{equation}
\left\Vert q\right\Vert_{_{L^{\infty}(\bar{I})}}^{2}\leq C\,(\ln
W)\left\Vert q\right\Vert_{_{1/2,I}}^{2}.\label{2.4}
\end{equation}
Here $C$ is a constant. Hence using (7.3) and (7.4) we obtain
\begin{align}
&\sum_{k=1}^{p}\sum_{i=1}^{I_{k}}\sum_{j=2}^{M}(\rho\mu_{k}^{M+1-j})^{-2\lambda_{k}}\Vert
r_{i,j}^{k}(\nu_{k},\phi_{k})\Vert_{2,\hat{\Omega}_{i,j}^{k}}^{2}+\sum_{l=1}^{L}\Vert
r_{^{l}}^{_{^{p+1}}}(\xi,\eta)\Vert_{_{_{2,S}}}^{2}\nonumber\\[.5pc]
&\quad\, \leq K(\ln
W)\Bigg(\sum_{k=1}^{p}\sum_{\gamma_{s}\subseteq\Omega^{k}\cup
B_{\rho}^{k},\mu(\hat{\gamma}_{s})<\infty}
d(A_{k},\gamma_{s})^{-2\lambda_{k}}\nonumber\\[.5pc]
&\qquad\, \times (\Vert
[u^{k}]\Vert_{_{0,\hat{\gamma}_{s}}}^{2}+\Vert
[(u_{\nu_{k}}^{k})^{a}]\Vert_{_{1/2,\hat{\gamma}_{s}}}^{2}+\Vert
[(u_{\phi_{k}}^{k})^{a}]\Vert_{_{1/2,\hat{\gamma}_{s}}}^{2})\nonumber\\[.5pc]
&\qquad\, + \sum_{\gamma_{s}\subseteq\Omega^{^{p+1}}}(\Vert
[u^{p+1}]\Vert_{_{0,\gamma_{s}}}^{2}+\Vert
[(u_{x_{1}}^{p+1})^{a}]\Vert_{_{1/2,\gamma_{s}}}^{2}+\Vert
[(u_{x_{2}}^{p+1})^{a}]\Vert_{_{1/2,\gamma_{s}}}^{2})\Bigg)\nonumber
\end{align}
\begin{align}
&\qquad\,
+\varepsilon_{_{W}}\,\Bigg(\sum_{k=1}^{p}\sum_{i=1}^{I_{k}}\sum_{j=2}^{M}(\rho\mu_{k}^{M+1-j})^{-2\lambda_{k}}\Vert
u_{i,j}^{k}(\nu_{k},\phi_{k})-h_{k}\Vert_{_{2,\hat{\Omega}_{i,j}^{k}}}^{2}\nonumber\\[.5pc]
&\qquad\, +\sum_{l=1}^{L}\Vert
u_{l}^{p+1}(\xi,\eta)\Vert_{_{2,S}}^{2}\Bigg).\label{8.5}
\end{align}
Let
\begin{align*}
&y_{i,j}^{k}(\nu_{k},\phi_{k})=u_{i,j}^{k}(\nu_{k},\phi_{k})+r_{i,j}^{k}(\nu_{k},\phi_{k})\,\:\textrm{and}\\[.5pc]
&y_{l}^{p+1}(\xi,\eta)=u_{l}^{p+1}(\xi,\eta)+r_{l}^{p+1}(\xi,\eta).
\end{align*}

Now we define a correction
$\{\{s_{i,j}^{k}(\nu_{k},\phi_{k})\}_{i,j,k},
\{s_{l}^{_{^{p+1}}}(\xi,\eta)\}_{l}\} $ such that $s_{i,1}^{k}=0$
for all $i$ and $k$, $s_{i,j}^{k}\in$
$H^{2}(\hat{\Omega}_{i,j}^{k})$ for $2\leq j\leq M$ and all $i$
and $k$, $s_{l}^{p+1}\in H^{2}(S)$ for $l=1,\dots,L$ and $w=y+s\in
H_{_{\beta}}^{^{2,2}}(\Omega).$

\begin{figure}[b]
\hskip 4pc{\epsfxsize=6.8cm\epsfbox{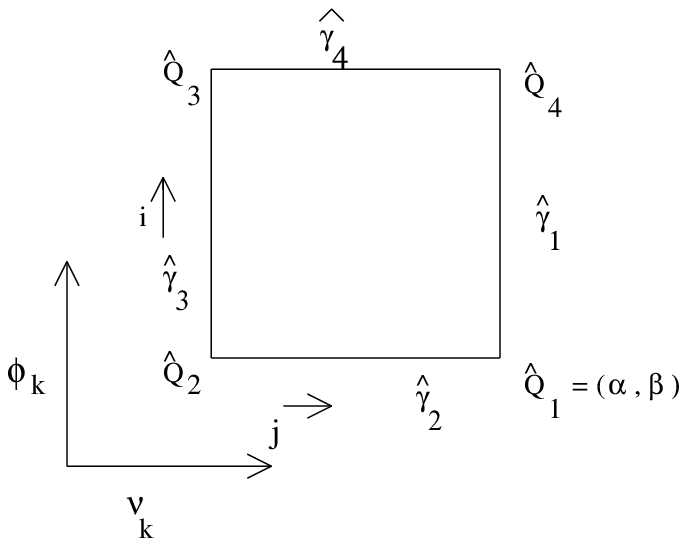}}\vspace{.7pc}
\end{figure}

Consider $\hat{\Omega}_{i,j}^{k}$ with $2\leq j<M.$ Let
\begin{align}
&F_{1}(\phi_{k})=-\frac{1}{2}(y_{i,j}^{k}-y_{i,j+1}^{k})\big|_{\hat{\gamma}_{1}},\nonumber\\[.5pc]
&G_{1}(\phi_{k})=-\frac{1}{2}(y_{i,j}^{k}-y_{i,j+1}^{k})_{\nu_{k}}\big|_{\hat{\gamma}_{1}},\,\,\textrm{and}\nonumber\\[.5pc]
&H_{1}(\phi_{k})=-\frac{1}{2}(y_{i,j}^{k}-y_{i,j+1}^{k})_{\phi_{k}}\big|_{\hat{\gamma}_{1}}.\label{2.6}
\end{align}
In the same way we define $F_{l},G_{l}$ and $H_{l}$ for
$l=1,\dots,4$. If $\gamma_{l}\,\, \subseteq \,\,\partial \Omega$
for some $l$, $F_{l},G_{l}$ and $H_{l}$ are defined to be
identically zero on $\hat{\gamma}_{l}.$ Now $F_{l},\, G_{l}$ and
$H_{l}$ are polynomials of degree $W$ that vanish at the end
points $\hat{Q}_{l}$ and $\,\hat{Q}_{l+1}$ of $\hat{\gamma}_{l}$.
If $\gamma_{3}\subseteq\partial\Omega_{i,1}^{k}
\cap\partial\Omega_{i,2}^{k}$ for some $i,k$ then the factor of
$\frac{1}{2}$ will be missing in the definition of
$F_{3}(\phi_{k}),\, G_{3}(\phi_{k})$ and $H_{3}(\phi_{k})$. We
wish to define $s_{i,j}^{k}(\nu_{k},\phi_{k})$ on
$\hat{\Omega}_{i,j}^{k}$ such that
$s_{i,j}^{k}\big|_{\hat{\gamma}_{l}}=F_{l},$
$(s_{i,j}^{k})_{\nu_{k}}\big|_{\hat{\gamma}_{l}}=G_{l}$ and
$(s_{i,j}^{k})_{\phi_{k}}\big|_{\hat{\gamma}_{l}}=H_{l}$ for
$l=1,\dots,4.$

We now cite Theorem~1.5.2.4 of \cite{grisvard}. The mapping
$u\rightarrow\{\{f_{k}\}_{k=0}^{m-1},\{g_{k}\}_{k=0}^{m-1}\}$
defined by $f_{k}= D_{\xi}^{k}u\big|_{\xi=0}$,
$g_{l}=D_{\eta}^{l}u\big|_{\eta=0}$ for $u\in
D(\overline{\mathbb{R}_{+}\times\mathbb{R}_{+}})$ has a unique
continuous extension as an operator from
$W_{p}^{m}(\overline{\mathbb{R}_{+}\times\mathbb{R}_{+}})$ onto
the subspace of
\begin{equation*}
T=\prod_{k=0}^{m-1}W^{m-k-1/p}(\mathbb{R}_{+})\times
\prod_{l=0}^{m-1}W^{m-l-1/p}(\mathbb{R}_{+})
\end{equation*}
defined by
\begin{enumerate}
\renewcommand\labelenumi{(\alph{enumi})}
\leftskip .2pc
\item $D_{\eta}^{l}f_{k}(0)=D_{\xi}^{k}g_{l}(0)$, $l+k<m-2/p$ for
all $p\neq2$, and

\item
$\int_{0}^{\delta}|D_{\eta}^{l}f_{k}(t)-D_{\xi}^{k}g_{l}(t)|^{2}{\rm
d}t/t<\infty$, $l+k=m-1$ for $p=2$.
\end{enumerate}

Hence using a partition of unity argument it is enough to show
that
\begin{enumerate}
\renewcommand\labelenumi{(\roman{enumi})}
\leftskip .2pc
\item $\int_{0}^{\delta}|D_{\phi_{k}}F_{1}(t+\beta)-
H_{2}(\alpha-t)|^{2}{\rm d}t/t$, and

\item $\int_{0}^{\delta}|G_{1}(t+\beta)-D_{\nu_{k}}F_{2}
(\alpha-t)|^{2}{\rm d}t/t$, are finite.
\end{enumerate}

Conditions (i) and (ii) follow by applying the above theorem to a
neighbourhood of the vertex $\hat{Q}_{1}=(\alpha,\beta)$ of
$\hat{\Omega}_{i,j}^{k}$.

Now
\begin{align*}
&\int_{0}^{\delta}|D_{\phi_{k}}F_{1}(t+\beta)-H_{2}(\alpha-t)|^{2}{\rm
d}t/t\\
&\quad\,
\leq2\int_{0}^{\delta}|D_{\phi_{k}}F_{1}(t+\beta)|^{2}{\rm
d}t/t+2\int_{0}^{\delta}|H_{2}(\alpha-t)|^{2}{\rm d}t/t\,.
\end{align*}

Moreover from Theorem~4.82 in \cite{schwab} we have that if $q(y)$
is a polynomial of degree $W$ on $I=(-1,1)$ such that
$q(-1)=q(1)=0$, then
\begin{align*}
\int_{-1}^{1}\frac{q^{2}(y)}{1-y^{2}}{\rm d}y\leq C\ln W\Vert
q\Vert_{_{L^{\infty}(I)}}^{2}.
\end{align*}
Now by (7.4),
\begin{align*}
\Vert q\Vert_{_{L^{\infty}(I)}}^{2}\leq K\ln W\Vert
q\Vert_{_{H^{1/2}(I)}}^{2}.
\end{align*}
Hence we conclude that
\begin{align}
\hskip -4pc{\rm (i)}\hskip 1pc
\int_{0}^{\delta}|D_{\phi_{k}}F_{1}(t+\beta)-H_{2}(\alpha-t)|^{2}{\rm
d}t/t\leq C(\ln W)^{2}(\Vert D_{\phi_{k}}F_{1}
\Vert_{_{1/2,\hat{\gamma}_{1}}}^{2}+\Vert H_{2}\Vert_{_{1/2,
\hat{\gamma}_{1}}}^{2}).\nonumber\\
\phantom{0}\hskip -1pc\label{2.7}
\end{align}
A similar result holds for (ii).

Hence we can define $\{\{s_{i,j}^{k}(\nu_{k},\phi_{k})\}_{i,j,k},
\{s_{l}^{_{^{p+1}}}(\xi,\eta)\}_{l}\}$ such that $s_{i,1}^{k} = 0$
for all $i$ and $k$, $s_{i,j}^{k}\in
H^{2}(\hat{\Omega}_{i,j}^{k})$ for $j\geq2$, $s_{l}^{p+1}\in
H^{2}(S)$ and $w=y+s\in H_{_{\beta}}^{^{2,2}}(\Omega)$.

Let $v_{i,j}^{k}(\nu_{k},\phi_{k})=r_{i,j}^{k}(\nu_{k},
\phi_{k})+s_{i,j}^{k}(\nu_{k},\phi_{k})$ and
$v_{l}^{p+1}(\xi,\eta)=r_{l}^{p+1}(\xi,\eta)+s_{l}^{p+1}(\xi,\eta)$.

Now from (7.7) we conclude that there is a constant $K$ such that
\begin{align}
&\sum_{k=1}^{p}\sum_{i=1}^{I_{k}}\sum_{j=2}^{M}(\rho\mu_{k}^{M+1-j})^{-2\lambda_{k}}\Vert
s_{i,j}^{k}(\nu_{k},\phi_{k})\Vert_{_{2,\hat{\Omega}_{i,j}^{k}}}^{2}+\sum_{l=1}^{L}\Vert
s_{l}^{p+1}(\xi,\eta)\Vert_{_{2,S}}^{2}\nonumber\\[.3pc]
&\quad\, \leq K(\ln
W)^{2}\left(\sum_{k=1}^{p}\sum_{\gamma_{s}\subseteq\Omega^{k}\cup
B_{\rho}^{k},\mu(\hat{\gamma}_{s})<\infty}
d(A_{k},\gamma_{s})^{-2\lambda_{k}} \right.\nonumber\\[.5pc]
&\qquad\, \times (\Vert
[u^{k}]\Vert_{_{0,\hat{\gamma}_{s}}}^{2}+\Vert
[(u_{\nu_{k}}^{k})^{a}]\Vert_{_{1/2,\hat{\gamma}_{s}}}^{2}+\Vert
[(u_{\phi_{k}}^{k})^{a}]\Vert_{_{1/2,\hat{\gamma}_{s}}}^{2})\nonumber\\[.3pc]
&\qquad\, \left. + \sum_{\gamma_{s}\subseteq\Omega^{^{p+1}}}(\Vert
[u^{p+1}]\Vert_{_{0,\gamma_{s}}}^{2}+\Vert
[(u_{x_{1}}^{p+1})^{a}]\Vert_{_{1/2,\gamma_{s}}}^{2}+\Vert
[(u_{x_{2}}^{p+1})^{a}]\Vert_{_{1/2,\gamma_{s}}}^{2})\right)\nonumber\\[.3pc]
&\qquad\,
+\varepsilon_{_{W}}\left(\sum_{k=1}^{p}\sum_{i=1}^{I_{k}}\sum_{j=2}^{M}(\rho\mu_{k}^{M+1-j})^{-2\lambda_{k}}\Vert
u_{i,j}^{k}(\nu_{k},\phi_{k})-h_{k}\Vert_{_{2,\hat{\Omega}_{i,j}^{k}}}^{2}
\right.\nonumber\\[.3pc]
&\qquad\,\left. +\sum_{l=1}^{L}\Vert
u_{l}^{p+1}(\xi,\eta)\Vert_{_{2,S}}^{2}\right).\label{1.8}
\end{align}
Combining (7.5) and (7.8) gives the estimate (7.1).\hfill
$\blacksquare$\vspace{.7pc}

We now prove the last result of this section.

\begin{lem}
Let $w=u+v\in H_{_{\beta}}^{^{2,2}}(\Omega).$ Here
$\{\{u_{i,j}^{k}(\nu_{k},\phi_{k})\}_{i,j,k},\{u_{l}^{_{^{p+1}}}(\xi,\eta)\}_{l}\}
\in\Pi^{^{M,W}}$ and
$\{\{v_{i,j}^{k}(\nu_{k},\phi_{k})\}_{i,j,k},\{v_{l}^{_{^{p+1}}}(\xi,\eta)\}_{l}\}$
is as defined in Lemma~$7.1$. Then the estimate
\begin{align}
&\Vert
w\Vert_{H_{_{\beta}}^{^{\frac{3}{2},\frac{3}{2}}}(\Gamma^{[0]})}^{2}+\left\Vert
\left(\frac{\partial w}{\partial
N}\right)_{A}\right\Vert_{H_{_{\beta}}^{^{\frac{1}{2},\frac{1}{2}}}(\Gamma^{[1]})}^{2}\nonumber\\[.5pc]
&\quad\, \leq C\,(\ln\, W)^{2}\left(\sum_{_{k\hbox{\rm :}\
\partial\Omega^{k}\cap\Gamma^{[0]}\neq\emptyset}}^{p}
|h_{k}|^{2}+\sum_{l\in\mathcal{D}}\sum_{k=l-1}^{l}\sum_{\gamma_{s}\subseteq\partial\Omega^{k}\cap\Gamma_{l},\mu(\hat{\gamma}_{s})<\infty}\right.\nonumber\\[.5pc]
&\qquad\,\times d(A_{k},\gamma_{s})^{-2\lambda_{k}}(\Vert
u^{k}-h_{k}\Vert_{0,\hat{\gamma}_{s}}^{2}+\Vert
u_{\nu_{k}}^{k}\Vert_{1/2,\hat{\gamma}_{s}}^{2}) \nonumber\\[.5pc]
&\qquad\, +
\sum_{l\in\mathcal{N}}\sum_{k=l-1}^{l}\sum_{\gamma_{s}\subseteq\partial\Omega^{k}\cap\Gamma_{l},\mu(\hat{\gamma}_{s})<\infty}
d(A_{k},\gamma_{s})^{-2\lambda_{k}}\left\Vert \left(\frac{\partial
u^{k}}{\partial
n}\right)_{\tilde{A}^{k}}^{a}\right\Vert_{1/2,\hat{\gamma}_{s}}^{2}\nonumber\\[.5pc]
 &\qquad\, +
\sum_{l\in\mathcal{D}}\sum_{\gamma_{s}\subseteq\partial\Omega^{^{p+1}}\cap\Gamma_{l}}
\left(\Vert u^{p+1}\Vert_{_{0,\gamma_{s}}}^{2} +\left\Vert
\left(\frac{\partial u^{p+1}}{\partial
T}\right)^{a}\right\Vert_{_{1/2,\gamma_{s}}}^{2}\right)\nonumber\\[.5pc]
&\qquad\,+\sum_{l\in\mathcal{N}}\sum_{\gamma_{s}\subseteq\partial\Omega^{^{p+1}}\cap\Gamma_{l}}\left\Vert
\left(\frac{\partial u^{p+1}}{\partial
N}\right)_{A}^{a}\right\Vert_{_{1/2,\gamma_{s}}}^{2}+
\sum_{k=1}^{p}\sum_{\gamma_{s}\subseteq\Omega^{k}\cup
B_{\rho}^{k},\mu(\hat{\gamma}_{s})<\infty}\nonumber
\end{align}
\begin{align}
&\quad\, \times d(A_{k},\gamma_{s})^{-2\lambda_{k}}(\Vert
[u^{k}]\Vert_{_{0,\hat{\gamma}_{s}}}^{2} +\Vert
[(u_{\nu_{k}}^{k})^{a}]\Vert_{_{1/2,\hat{\gamma}_{s}}}^{2}+\Vert
[(u_{\phi_{k}}^{k})^{a}]\Vert_{_{1/2,\hat{\gamma}_{s}}}^{2})\nonumber\\[.5pc]
&\quad\, \left. + \sum_{\gamma_{s}\subseteq\Omega^{^{p+1}}}(\Vert
[u^{p+1}]\Vert_{_{0,\gamma_{s}}}^{2}+\Vert
[(u_{x_{1}}^{p+1})^{a}]\Vert_{_{1/2,\gamma_{s}}}^{2}+\Vert
[(u_{x_{2}}^{p+1})^{a}]\Vert_{_{1/2,\gamma_{s}}}^{2})\right)\nonumber\\[.5pc]
&\quad\, +
\varepsilon_{_{W}}\,\left(\sum_{k=1}^{p}|h_{k}|^{2}+\sum_{k=1}^{p}\sum_{i=1}^{I_{k}}\sum_{j=2}^{M}(\rho\mu_{k}^{M+1-j})^{-2\lambda_{k}}\Vert
u_{i,j}^{k}(\nu_{k},\phi_{k})-h_{k}\Vert_{_{2,\hat{\Omega}_{i,j}^{k}}}^{2}\right.\nonumber\\[.5pc]
&\quad\,\left.+\sum_{l=1}^{L}\Vert
u_{l}^{p+1}(\xi,\eta)\Vert_{_{2,S}}^{2}\right)\label{8.9}
\end{align}
holds. Here $\varepsilon_{_{W}}$ is exponentially small in $W.$
Now
\begin{equation*}
\Vert
w\Vert_{_{H_{_{\beta}}^{^{\frac{3}{2},\frac{3}{2}}}(\Gamma^{[0]})}}=\inf_{q|_{\Gamma^{[0]}}=w}\{\Vert
q\Vert_{_{H_{_{\beta}}^{^{2,2}}(\Omega)}}\}.
\end{equation*}
\end{lem}
Let $\theta_{k}\in C^{2}(\mathbb{R})$ such that $\theta_{k}=1$ for
$r_{k}\leq\rho\mu_{k}$ and $\theta_{k}=0$ for $r_{k}\geq\rho.$ Let
$q_{k}=q\theta_{k}$ and $q_{0}=1-{\sum_{k=1}^{p}q_{k}}$. Let
$\theta_{0}=1-{\sum_{k=1}^{p}\theta_{k}}.$ Define
$\Omega_{_{\rho\mu_{k}}}^{^{k}}=\{x\hbox{\rm :}\
d(A_{k},x)<\rho\mu_{k}\} $ for $k=1,\dots,p$ and let
$\tilde{\Omega}^{p+1}=\Omega\setminus{\bigcup_{k=1}^{p}\bar{\Omega}_{_{\rho\mu_{k}}}^{k}}.$
Then it can be concluded that
\begin{align}
&\hskip -4pc \Vert
w\Vert_{_{H_{_{\beta}}^{\frac{3}{2},\frac{3}{2}}(\Gamma^{[0]})}}^{2}\nonumber\\[.5pc]
&\hskip -4pc \quad\, \leq C\left(\sum_{k\hbox{\rm :}\
\partial\Omega^{k}\cap\Gamma^{[0]}\neq\emptyset} \
\inf_{_{_{q_{k}|_{\partial\Omega^{k}\cap\Gamma^{[0]}}=w\theta_{k}}}}\{\Vert
q_{k}\Vert_{_{H_{_{\beta}}^{^{2,2}}(\Omega^{k})}}^{2}\}
+\inf_{_{_{q_{_{0}}|_{\partial\tilde{\Omega}^{^{p+1}}\cap\Gamma^{[0]}}=w\theta_{0}}}}\{\Vert
q_{0}\Vert_{_{H^{2}(\tilde{\Omega}^{^{p+1}})}}^{2}\}
\right).\label{2.10}
\end{align}
Now using (3.9) we have
\begin{equation}
\Vert q_{k}\Vert_{_{H_{_{\beta}}^{^{2,2}}(\Omega^{k})}}^{2}\leq
C(|h_{k}|^{2}+\Vert
(q_{k}(\nu_{k},\phi_{k})-h_{k})\hbox{e}^{-2(1-\beta_{k})\nu_{k}}\Vert
_{_{2,\hat{\Omega}^{k}}}^{2}).\label{2.11}
\end{equation}
Let us choose the cut-off function $\theta_{k}$ to be a piecewise
polynomial such that
\begin{align*}
\theta_{k}(\nu_{k})&=1\quad\hbox{for} \ \ \nu_{k}\leq
\ln(\rho\mu_{k}),\\[.5pc]
\theta_{k}(\ln(\rho\mu_{k}))&=1,\theta_{k}^{(1)}(\ln(\rho\mu_{k}))=0,\theta_{k}^{(2)}
(\ln(\rho\mu_{k}))=0,\\[.5pc]
\theta_{k}(\ln\rho)&=0,\theta_{k}^{(1)}(\ln\rho)=0,\theta_{k}^{(2)}(\ln\rho)=0,\quad\hbox{and}\\[.5pc]
\theta_{k}(\nu_{k})&=0\quad\hbox{for} \ \ \nu_{k}\geq \ln\rho.
\end{align*}
Here $\theta_{k}^{(l)}$ denotes the $l$th derivative of
$\theta_{k}$ with respect to $\nu_{k}.$ Then $\theta_{k}$ is a
polynomial of degree five in $\nu_{k}$ for
$\ln(\rho\mu_{k})\leq\nu_{k}\leq \ln\rho$. Now using (7.10) and
(7.11) we have\pagebreak
\begin{align*}
\hskip -4pc&\sum_{k\hbox{\rm :}\
\partial\Omega^{k}\cap\Gamma^{[0]}\neq\emptyset}\
\inf_{_{q_{k}|_{\partial\Omega^{k}\cap
\Gamma^{[0]}}=w\theta_{k}}}\{\Vert
q_{k}\Vert_{_{H_{_{\beta}}^{^{2,2}}(\Omega^{k})}}^{2}\}\\[.5pc]
\hskip -4pc&\quad\,\leq C\left(\left(\sum_{k\hbox{\rm :}\
\partial\Omega^{k}\cap\Gamma^{[0]}}|h_{k}|^{2}+\sum_{k\hbox{\rm :}\ \Gamma_{k}\cap\Gamma^{[0]}\neq\emptyset}
\Vert (q_{k}(\nu_{k},\psi_{l}^{k})-h_{k}){\rm e}^{-2(1-\beta_{k})\nu_{k}}\Vert_{_{3/2,(-\infty,ln\rho)}}^{2}\right)\right.\\[.5pc]
\hskip -4pc&\qquad\,+\left.\sum_{k\hbox{\rm :}\
\Gamma_{k+1}\cap\Gamma^{[0]}\neq\emptyset}\Vert
(q_{k}(\nu_{k},\psi_{u}^{k})-h_{k}){\rm
e}^{-2(1-\beta_{k})\nu_{k}}\Vert_{_{3/2,(-\infty,ln\rho)}}^{2}\right).
\end{align*}
Let $\eta_{j}^{k}=\ln\,\rho+(M+1-j)\,\, \ln\,\mu_{k}$ and
$I_{j}^{k}=(\eta_{j-1}^{k},\eta_{j}^{k})$. Then
\begin{align}
&\Vert
(q_{_{k}}(\nu_{_{k}},\psi_{u}^{k})-h_{_{k}})\hbox{e}^{-2(1-\beta_{k})\nu_{k}}\Vert_{_{3/2,(-\infty,ln\rho)}}^{2}\nonumber\\[.5pc]
&\quad\,\leq C\left\{\sum_{j=2}^{M+1}\Vert
(q_{k}(\nu_{k},\psi_{u}^{k})-h_{_{k}})\hbox{e}^{^{-2(1-\beta_{k})\nu_{k}}}\Vert_{_{0,(\eta_{j-1}^{k},\eta_{j}^{k})}}^{^{2}}\right.\nonumber\\[.5pc]
&\qquad\,\left.+\sum_{j=2}^{M+1}\Vert ((q_{k}(\nu_{k},\psi_{u}^{k})-h_{_{k}})\hbox{e}^{^{-2(1-\beta_{k})\nu_{k}}})_{\nu_{k}}\Vert_{_{1/2,(\eta_{j-1}^{k},\eta_{j}^{k})}}^{^{2}}\right.\nonumber\\[.5pc]
&\qquad\,\left.+\sum_{j=2}^{M+1}\int_{0}^{\delta}\left|\frac{{\rm
d}}{{\rm
d}s}((q_{_{k}}(s,\psi_{u}^{k})-h_{_{k}})\hbox{e}^{^{-(1-\beta_{k})s}})(\eta_{j}^{k}+\sigma)\right.\right.\nonumber\\[.5pc]
&\qquad\,\left.\left.-\frac{{\rm d}}{{\rm
d}s}((q_{_{k}}(s,\psi_{u}^{k})-h_{_{k}})
\hbox{e}^{^{-(1-\beta_{k})s}})(\eta_{j}^{k}-\sigma)\right|^{2}\frac{{\rm
d}\sigma}{\sigma}\right\}.\label{8.12}
\end{align}
Here $\delta>0$. Now
\begin{align*}
&\int_{0}^{\delta}\left|\frac{{\rm d}}{{\rm
d}s}((q_{k}(s,\psi_{u}^{k})-h_{k})
\hbox{e}^{-(1-\beta_{k})s})(\eta_{j}^{k}+\sigma)\right.\\[.5pc]
&\qquad\,\left.-\frac{{\rm d}}{{\rm d}s}((q_{k}(s,\psi_{u}^{k})-h_{k}){\rm e}^{-(1-\beta_{k})s})(\eta_{j}^{k}-\sigma)\right|^{2}\frac{{\rm d}\sigma}{\sigma}\\[.5pc]
&\quad\,\leq
K\left(\sum_{l=0}^{1}\int_{0}^{\delta}\left|\frac{{\rm
d}^{l}}{{\rm d}s^{l}}(q_{k}-h_{k})
(\eta_{j}^{k}+\sigma,\psi_{u}^{k})\right.\right.\\[.5pc]
&\qquad\,\left.\left.-\frac{{\rm d}^{l}}{{\rm d}s^{l}}(q_{k}-h_{k})(\eta_{j}^{k}-\sigma,\psi_{u}^{k})\right|^{2}\hbox{e}^{-2(1-\beta_{k})(\eta_{j}^{k}+\sigma)}\frac{{\rm d}\sigma}{\sigma}\right.\\[.5pc]
&\qquad\,+ \left.\sum_{l=0}^{1}\int_{0}^{\delta}\left|\frac{{\rm
d}^{l}}{{\rm
d}s^{l}}(q_{k}-h_{k})(\eta_{j}^{k}-\sigma,\psi_{u}^{k})\right|^{2}\right.\\[.5pc]
&\qquad\,\left.\times({\rm
e}^{-2(1-\beta_{k})(\eta_{j}^{k}+\sigma)}-{\rm
e}^{-2(1-\beta_{k})(\eta_{j}^{k}-\sigma)})\frac{{\rm
d}\sigma}{\sigma}\right).
\end{align*}
Hence
\begin{align}
&\Vert (q_{k}(\nu_{k},\psi_{u}^{k})-h_{k}){\rm
e}^{-2(1-\beta_{k})\nu_{k}}
\Vert_{_{3/2,(-\infty,\ln\rho)}}^{2}\nonumber\\[.3pc]
&\quad\,\leq
C\left(\sum_{\gamma_{s}\subseteq\Gamma^{[0]}\cap\Gamma_{k+1},\mu(\hat{\gamma}_{s})<\infty}
d(A_{k},\gamma_{s})^{-2\lambda_{k}}\right.\nonumber\\[.3pc]
&\qquad\,\left.\times\left(\Vert
(q_{k}-h_{k})\Vert_{_{0,\hat{\gamma}_{s}}}^{2}+\left\Vert
\frac{{\rm d}q_{k}}{{\rm
d}\nu_{k}}\right\Vert_{_{1/2,\hat{\gamma}_{s}}}^{2}+\ln\, W\Vert
(q_{k}-h_{k})\Vert_{_{1,\infty,\hat{\gamma}_{s}}}^{2}\right)\right).\label{2.13}
\end{align}
This follows from Theorem 4.82 in \cite{schwab} which states that
if $p(y)$ is a polynomial of degree $N$ in $y$ such that
$p(1)=p(-1)=0,$ then
\begin{equation*}
\int_{-1}^{1}\frac{p^{2}(y)}{1-y^{2}}{\rm d}y\leq\, C\,\, \ln
N\Vert p\Vert_{_{L^{\infty}(\bar{I})}}^{2}.
\end{equation*}
Now $q_{k}(\nu_{k},\psi_{u}^{k})-h_{k}=\theta_{k}\,
w(\nu_{k},\psi_{u}^{k})-h_{k}=\theta_{k}(w(\nu_{k},\psi_{u}^{k})-h_{k})+(\theta_{k}-1)h_{k}.$

Moreover $w=r+s+u$, as has been defined in Lemma~7.1. Here
$w(\nu_{k},\psi_{u}^{k})$ is a polynomial of degree $W,$
$s(\nu_{k},\psi_{u}^{k})=0$ and $r(\nu_{k},\psi_{u}^{k})$ is a
polynomial degree four.

Hence using (7.4) and (7.13) we conclude that
\begin{align}
&\Vert(q_{k}(\nu_{k},\psi_{u}^{k})-h_{k}){\rm
e}^{-2(1-\beta_{k})\nu_{k}}\Vert_{_{3/2,(-\infty,\ln\rho)}}^{2}\nonumber\\[.5pc]
&\quad\,\leq C\left(
\sum_{\gamma_{s}\subseteq\Gamma_{k+1}\bigcap\partial\Omega^{k}}
d(A_{k},\gamma_{s})^{-2\lambda_{k}} ((\ln W)^{2}(\Vert
(u(\nu_{k},\psi_{u}^{k})-h_{k})\Vert_{_{3/2,\hat{\gamma}_{s}}}^{2}
\right.\nonumber\\[.3pc]
&\qquad\,\left. +|h_{k}|^{2})+\ln W\Vert
r(\nu_{k},\psi_{u}^{k})\Vert_{_{3/2,\hat{\gamma}_{s}}}^{2})\right).\label{2.14}
\end{align}
Hence using (7.5) and (7.14) it can be concluded that
\begin{align}
&\sum_{k\hbox{\rm :}\
\partial\Omega^{k}\cap\Gamma^{[0]}\neq\emptyset} \
\inf_{_{q_{k}|_{\partial\Omega^{k}\cap\Gamma^{[0]}}=w\theta_{k}}}\{\Vert
q_{k}\Vert_{_{H_{_{\beta}}^{^{2,2}}(\Omega^{k})}}^{2}\}\nonumber\\[.5pc]
&\quad\,\leq C(\ln
W)^{2}\left(\sum_{k=1}^{p}|h_{k}|^{2}+\sum_{l\in\mathcal{D}}\sum_{k=l-1}^{l}\sum_{\gamma_{s}\subseteq\partial\Omega^{k}\cap\Gamma_{l},\mu(\hat{\gamma}_{s})<\infty}
d(A_{k},\gamma_{s})^{-2\lambda_{k}}\right.\nonumber\\[.5pc]
&\qquad\,\left.\times(\Vert
u^{k}-h_{k}\Vert_{_{0,\hat{\gamma}_{s}}}^{2}+\Vert
u_{\nu_{k}}^{k}\Vert_{_{1/2,\hat{\gamma}_{s}}}^{2})+
\sum_{k=1}^{p}\sum_{\gamma_{s}\subseteq\Omega^{k}\cup
B_{\rho}^{k},\mu(\hat{\gamma}_{s})<\infty}
d(A_{k},\gamma_{s})^{-2\lambda_{k}} \right.\nonumber\\[.5pc]
&\qquad\, \times(\Vert
[u^{k}]\Vert_{_{0,\hat{\gamma}_{s}}}^{2}+\Vert
[(u_{\nu_{k}}^{k})^{a}]\Vert_{_{1/2,\hat{\gamma}_{s}}}^{2}+\Vert
[(u_{\phi_{k}}^{k})^{a}]\Vert_{_{1/2,\hat{\gamma}_{s}}}^{2})\nonumber\\[.3pc]
&\qquad\,\left. + \sum_{\gamma_{s}\subseteq\Omega^{^{p+1}}}(\Vert
[u^{p+1}] \Vert_{_{0,\gamma_{s}}}^{2}+\Vert
[(u_{x_{1}}^{p+1})^{a}]\Vert_{_{1/2,\gamma_{s}}}^{2}+\Vert
[(u_{x_{2}}^{p+1})^{a}]
\Vert_{_{1/2,\gamma_{s}}}^{2})\right)\nonumber 
\end{align}
\begin{align}
&\qquad\, +
\varepsilon_{_{W}}\left(\sum_{k=1}^{p}\sum_{j=2}^{M}\sum_{i=1}^{I_{k}}(\rho\mu_{k}^{M+1-j})^{-2\lambda_{k}}\Vert
u_{i,j}^{k}(\nu_{k},\phi_{k})-h_{k}\Vert_{_{2,\hat{\Omega}_{i,j}^{k}}}^{2}\right.\nonumber\\[.5pc]
&\qquad\,\left.+\sum_{k=1}^{p}|h_{k}|^{2}+\sum_{l=1}^{L}\Vert
u_{l}^{p+1}(\xi,\eta)\Vert_{_{2,S}}^{2}\right).\label{2.15}
\end{align}
In the same way we can show that
\begin{align}
&\inf_{_{q_{_{0}}|_{\partial\tilde{\Omega}^{p+1}\cap\Gamma^{[0]}}=w\theta_{0}}}(\Vert
q_{0}\Vert_{_{H^{2}(\tilde{\Omega}^{p+1})}}^{2})\nonumber\\[.5pc]
&\quad\, \leq  C(\ln W)^{2}\left(\sum_{_{k\hbox{\rm :}\
\partial\Omega^{k}\cap\Gamma^{[0]}\neq\emptyset}}^{p} |h_{k}|^{2}
+\sum_{l\in\mathcal{D}}\sum_{k=l-1}^{l}\sum_{\gamma_{s}\subseteq\partial\Omega^{k}
\cap\Gamma_{l},\mu(\hat{\gamma}_{s})<\infty}\right.\nonumber\\[.5pc]
&\qquad\,\left. \times d(A_{k},\gamma_{s})^{-2\lambda_{k}} (\Vert
u^{k}-h_{k}\Vert_{0,\hat{\gamma}_{s}}^{2}  +\Vert
u_{\nu_{k}}^{k}\Vert_{1/2,\hat{\gamma}_{s}}^{2})+
\sum_{k=1}^{p}\sum_{\gamma_{s}\subseteq\Omega^{k}\cup
B_{\rho}^{k},\mu(\hat{\gamma}_{s}) <\infty}\right.\nonumber\\[.5pc]
&\qquad\,\times d(A_{k},\gamma_{s})^{-2\lambda_{k}}(\Vert
[u^{k}]\Vert_{_{0,\hat{\gamma}_{s}}}^{2}+\Vert
[(u_{\nu_{k}}^{k})^{a}]\Vert_{_{1/2,\hat{\gamma}_{s}}}^{2}+\Vert
[(u_{\phi_{k}}^{k})^{a}]\Vert_{_{1/2,\hat{\gamma}_{s}}}^{2})\nonumber\\[.5pc]
&\qquad\, +\sum_{\gamma_{s}\subseteq\Omega^{^{p+1}}}(\Vert
[u^{p+1}]\Vert_{_{0,\gamma_{s}}}^{2}+\Vert
[(u_{x_{1}}^{p+1})^{a}]\Vert_{_{1/2,\gamma_{s}}}^{2}
+\Vert [(u_{x_{2}}^{p+1})^{a}]\Vert_{_{1/2,\gamma_{s}}}^{2})\nonumber\\[.5pc]
&\qquad\,\left. +
\sum_{l\in\mathcal{D}}\sum_{\gamma_{s}\subseteq\partial\Omega^{^{p+1}}\cap\Gamma_{l}}
\left(\Vert u^{p+1}\Vert_{_{0,\gamma_{s}}}^{2}+\left\Vert
\left(\frac{\partial u^{p+1}}{\partial T}
\right)^{a}\right\Vert_{_{1/2,\gamma_{s}}}^{2}\right)\right)\nonumber
\\[.5pc]
&\qquad\, +
\varepsilon_{_{W}}\left(\sum_{k=1}^{p}\sum_{j=2}^{M}\sum_{i=1}^{I_{k}}(\rho\mu_{k}^{M+1-j})^{-2\lambda_{k}}\Vert
u_{i,j}^{k}(\nu_{k},\phi_{k})-h_{k}\Vert_{_{2,\hat{\Omega}_{i,j}^{k}}}^{2}\right.\nonumber\\[.5pc]
&\qquad\,\left.+\sum_{k=1}^{p}|h_{k}|^{2}+\sum_{l=1}^{L}\Vert
u_{l}^{p+1}(\xi,\eta)\Vert_{_{2,S}}^{2}\right).\label{2.16}
\end{align}
Now
\begin{equation*}
\left\Vert \left(\frac{\partial w}{\partial
N}\right)_{A}\right\Vert_{_{H_{\beta}^{\frac{1}{2},\frac{1}{2}}(\Gamma^{[1]})}}^{2}
=\inf_{_{q|_{\Gamma^{[1]}}=\big(\frac{\partial w}{\partial
N}\big)_{A}}}(\left\Vert
q\right\Vert_{H_{\beta}^{^{1,1}}(\Omega)}^{2}).
\end{equation*}
Here
\begin{equation*}
\Vert q\Vert_{_{H_{_{\beta}}^{^{1,1}}(\Omega)}}^{2}=\Vert
q\Vert_{_{L^{2}(\Omega)}}^{2}+\sum_{|\alpha|=1}\Vert
\Phi_{\beta}D^{\alpha}q\Vert_{_{L^{2}(\Omega)}}^{2}.
\end{equation*}
Let $\theta_{k}\in C^{2}(\mathbb{R})$ be as defined earlier and
$q_{k}=q\,\theta_{k}.$ Let $q_{0}=1-{ \sum_{k=1}^{p}q_{k}}.$ Then,
as before
\begin{align*}
&\hskip -4pc \left\Vert \left(\frac{\partial w}{\partial
N}\right)_{A}\right\Vert_{_{H_{_{\beta}}^{^{\frac{1}{2},\frac{1}{2}}}(\Gamma^{[1]})}}^2\\[.5pc]
&\hskip -4pc \quad\,\leq
C\left(\inf_{q_{k}|_{\partial\Omega^{k}\cap\Gamma^{[1]}}=\big(\frac{\partial
w}{\partial N}\big)_{_{A}}\theta_{k}}(\left\Vert
q_{k}\right\Vert_{_{H_{_{\beta}}^{^{1,1}}(\Omega^{k})}}^{2})+\inf_{q_{0}|_{\partial\tilde{\Omega}^{p+1}\cap\Gamma^{[1]}}=\big(\frac{\partial
w}{\partial N}\big)_{_{A}}\theta_{0}}(\left\Vert
q_{0}\right\Vert_{_{H^{1}(\tilde{\Omega}^{p+1})}}^{2})\right).
\end{align*}
Now
\begin{equation*}
\int_{\hat{\Omega}^{k}}\hbox{e}^{2\beta_{k}\nu_{k}}|q(\nu_{k},\phi_{k})|^{2}\hbox{d}\nu_{k}\hbox{d}\phi_{k}\leq
K\,\,\left\Vert q\right\Vert
_{_{H_{_{\beta}}^{^{1,1}}(\Omega^{k})}}^{2}
\end{equation*}
for $\beta_{k}>0.$ And so
\begin{equation*}
\int_{\hat{\Omega}^{k}}\hbox{e}^{2\beta_{k}\nu_{k}}|q_{k}(\nu_{k},\phi_{k})|^{2}\hbox{d}\nu_{k}\hbox{d}\phi_{k}\leq
K\,\left\Vert
q_{k}\right\Vert_{_{H_{_{\beta}}^{^{1,1}}(\Omega^{k})}}^{2}.
\end{equation*}
Hence for $0<\beta_{k}<1$, there exists a constant $C$ such that
\begin{align*}
&\frac{1}{C}\left(\sum_{k=1}^{p}\inf_{q_{k}|_{\partial\hat{\Omega}^{k}\cap\Gamma^{[1]}}={\rm
e}^{-\nu_{k}}\big(\frac{\partial w}{\partial
n}\big)_{\tilde{A}^{k}}\theta_{k}}\right.\\[.5pc]
&\qquad\,\left.\times\left(\int_{\hat{\Omega}^{k}}{\rm e}^{2\beta_{k}\nu_{k}}\sum_{|\alpha|\leq1}|D_{\nu_{k},\phi_{k}}^{\alpha}q_{k}(\nu_{k},\phi_{k})|^{2}\hbox{d}\nu_{k}\hbox{d}\phi_{k}\right)\right)\\[.5pc]
&\quad\,\leq
\sum_{k=1}^{p}\inf_{q_{k}|_{\partial\hat{\Omega}^{k}\cap\Gamma^{[1]}}={\rm
e}^{-\nu_{k}}\big(\frac{\partial w}{\partial
n}\big)_{\tilde{A}^{k}}\theta_{k}}(\Vert
q_{k}\Vert_{_{H_{_{\beta}}^{^{1,1}}(\Omega^{k})}}^{2})\\[.5pc]
&\quad\,\leq
C\left(\sum_{k=1}^{p}\inf_{q_{k}|_{\partial\hat{\Omega}^{k}\cap\Gamma^{[1]}}={\rm
e}^{-\nu_{k}}\big(\frac{\partial w}{\partial
n}\big)_{\tilde{A}^{k}}\theta_{k}}\right.\\[.5pc]
&\qquad\,\left.\times\left(\int_{\hat{\Omega}^{k}}\hbox{e}^{2\beta_{k}\nu_{k}}\sum_{|\alpha|\leq1}|D_{\nu_{k},\phi_{k}}^{\alpha}q_{k}(\nu_{k},\phi_{k})|^{2}\hbox{d}\nu_{k}\hbox{d}\phi_{k}\right)\right).
\end{align*}
Thus by similar arguments as before it can be shown that
\begin{align}
&\left\Vert \left(\frac{\partial w}{\partial
N}\right)_{A}\right\Vert_{_{H_{_{\beta}}^{^{\frac{1}{2},\frac{1}{2}}}(\Gamma^{[1]})}}^{2}\nonumber\\[.5pc]
&\quad\,\leq
C(\hbox{ln}W)^{2}\left(\sum_{l\in\mathcal{N}}\sum_{k=l-1}^{l}\sum_{\gamma_{s}\subseteq\partial\Omega^{k}\cap\Gamma_{l},\mu(\hat{\gamma}_{s})<\infty}
d(A_{k},\gamma_{s})^{-2\lambda_{k}}\right.\nonumber
\end{align}
\begin{align}
&\qquad\,\times\left\Vert \left(\frac{\partial u^{k}}{\partial
n}\right)_{\tilde{A}^{k}}^{a}\right\Vert_{_{1/2,\hat{\gamma}_{s}}}^{2}+ \sum_{l\in\mathcal{N}}\sum_{\gamma_{s}\subseteq\partial\Omega^{^{p+1}}\cap\Gamma_{l}}  \left\Vert \left(\frac{\partial u^{p+1}}{\partial N}\right)_{A}^{a}\right\Vert_{_{1/2,\gamma_{s}}}^{2}\nonumber\\[.5pc]
&\qquad\,+ \sum_{k=1}^{p}\sum_{\gamma_{s}\subseteq\Omega^{k}\cup
B_{\rho}^{k},\mu(\hat{\gamma}_{s})<\infty}
d(A_{k},\gamma_{s})^{-2\lambda_{k}}\nonumber\\[.5pc]
&\qquad\,\times(\Vert
[u^{k}]\Vert_{_{0,\hat{\gamma}_{s}}}^{2}+\Vert
[(u_{\nu_{k}}^{k})^{a}]\Vert_{_{1/2,\hat{\gamma}_{s}}}^{2}+\Vert
[(u_{\phi_{k}}^{k})^{a}]\Vert_{_{1/2,\hat{\gamma}_{s}}}^{2})\nonumber\\[.5pc]
&\qquad\,+ \left.\sum_{\gamma_{s}\subseteq\Omega^{^{p+1}}}(\Vert
[u^{p+1}]\Vert_{_{0,\gamma_{s}}}^{2}+\Vert
[(u_{x_{1}}^{p+1})^{a}]\Vert_{_{1/2,\gamma_{s}}}^{2}+\Vert
[(u_{x_{2}}^{p+1})^{a}]\Vert_{_{1/2,\gamma_{s}}}^{2})\right)\nonumber\\[.5pc]
&\qquad\,+
\varepsilon_{_{W}}\left(\sum_{k=1}^{p}\sum_{j=2}^{M}\sum_{i=1}^{I_{k}}(\rho\mu_{k}^{M+1-j})^{-2\lambda_{k}}\Vert
u_{i,j}^{k}(\nu_{k},\phi_{k})-h_{k}\Vert_{2,\hat{\Omega}_{i,j}^{k}}^{2}\right.\nonumber\\[.5pc]
&\qquad\,\left.+\sum_{k=1}^{p}|h_{k}|^{2}+\sum_{l=1}^{L}\Vert
u_{l}^{p+1}(\xi,\eta)\Vert_{_{2,S}}^{2}\right)
\end{align}
Combining (7.15)--(7.17) we obtain the required result.\hfill
$\blacksquare$

\section*{Acknowledgement}

This research is partly supported by CDAC (Center for Development
of Advanced Computing, Pune).


\begin{thebibliography}{99}

\bibitem{babcra} Babuska I, Craig A, Mandel J and Pitkaranta J, Efficient
preconditioning for the p version of the finite element method in
two dimensions, {\it SIAM J.~Num. Anal.} {\bf 28} (1991) 624

\bibitem{babguo1} Babuska I and Guo B Q, Regularity of the solution of elliptic
problems with piecewise analytic data, Part-I, {\it SIAM J.~Math.
Anal.} {\bf 19} (1988) 172--203

\bibitem{babguo2} Babuska I and Guo B Q, The h-p version of the finite element
method on domains with curved boundaries, {\it SIAM J.~Num. Anal.}
{\bf 25} (1988) 837--861

\bibitem{babguo3} Babuska I and Guo B Q, Regularity of the solution of elliptic
problems with piecewise analytic data, Part-II, {\it SIAM J.~Math.
Anal.} {\bf 20} (1989) 763--781

\bibitem{prsb} Dutt P K and Bedekar S, Spectral methods for hyperbolic initial
boundary value problems on parallel computers, {\it J.~Comput.
Appl. Math.} {\bf 134} (2001) 165--190

\bibitem{duttora1} Dutt P, Tomar S and Kumar R, Stability estimates for h-p
spectral element methods for elliptic problems, {\it Proc. Indian
Acad. Sci (Math. Sci.)} {\bf 112(4)} (2002) 601--639

\bibitem{duttom} Dutt P and Tomar S, Stability estimates for h-p spectral element
methods for general elliptic problems on curvilinear domains, {\it
Proc. Indian Acad. Sci (Math. Sci.)} {\bf 113} (2003) 395--429

\bibitem{grisvard} Grisvard P, Elliptic problems in non-smooth domains (Pitman
Advanced Publishing Program) (1985)

\bibitem{babguo4} Guo B Q and Babuska I, On the regularity of elasticity problems
with piecewise analytic data, {\it Adv. Appl. Math.} {\bf 14}
(1993) 307--347

\bibitem{guocao} Guo B and Cao W, A preconditioner for the h-p version of
the finite element method in two dimensions, {\it Num. Math.} {\bf
75} (1996) 59

\bibitem{karnia} Karniadakis G and Spencer Sherwin J, Spectral/hp element
methods for CFD, (Oxford University Press) (1999)

\bibitem{schwab} Schwab Ch, p and h-p Finite element methods (Oxford: Clarendon Press)
(1998)

\bibitem{duttora2} Tomar S K, Dutt P and Rathish Kumar B V, An efficient
and exponentially accurate parallel h-p spectral element method
for elliptic problems on polygonal domains--The Dirichlet case,
Lecture Notes in Computer Science 2552, High Performance Computing
HiPC (Springer-Verlag) (2002)

\bibitem{tomarth} Tomar S K, h-p Spectral element methods for elliptic problems on
non-smooth domains using parallel computers, Ph.D. thesis (India:
IIT Kanpur) (2001); Reprint available as Tec.~Rep.~no.~1631,
Department of Applied Mathematics, University of Twente, The
Netherlands. http://www.math.utwente.nl/publications
\end{thebibliography}
\end{document}